\documentclass[reqno,11pt, psfig]{amsart}
\usepackage{amssymb, amsmath,amssymb, amscd, amsbsy,latexsym}
\usepackage{amsthm, amsfonts,mathrsfs}
\usepackage[T1]{fontenc}
\usepackage{times}
\usepackage{epsfig}
\usepackage{color}
\usepackage[all]{xypic}
\input amssym.def
\input amssym.tex

\def\eqdefa{\buildrel\hbox{\footnotesize def}\over =}

\def\inte#1{
\displaystyle\mathop{#1\kern0pt}^\circ }

\let\pa=\partial

\let\D=\Delta

\def\cB{{\mathcal B}}
\def\cC{{\mathcal C}}

\def\cF{{\mathcal F}}

\def\pa{\partial}

\def\virgp{\raise 2pt\hbox{,}}
\def\cdotpv{\raise 2pt\hbox{;}}

\def\eqdefa{\buildrel\hbox{\footnotesize def}\over =}

\def\C{\mathop{\mathbb C\kern 0pt}\nolimits}
\def\DD{\mathop{\mathbb D\kern 0pt}\nolimits}
\def\EE{\mathop{{\mathbb E \kern 0pt}}\nolimits}
\def\K{\mathop{\mathbb K\kern 0pt}\nolimits}
\def\N{\mathop{\mathbb N\kern 0pt}\nolimits}
\def\Q{\mathop{\mathbb Q\kern 0pt}\nolimits}
\def\R{\mathop{\mathbb R\kern 0pt}\nolimits}
\def\SS{\mathop{\mathbb S\kern 0pt}\nolimits}
\def\ZZ{\mathop{\mathbb Z\kern 0pt}\nolimits}
\def\TT{\mathop{\mathbb T\kern 0pt}\nolimits}
\def\P{\mathop{\mathbb P\kern 0pt}\nolimits}

\def\Supp{\mathop{\rm Supp}\nolimits\ }

\newcommand{\beq}{\begin{equation}}
\newcommand{\eeq}{\end{equation}}
\newcommand{\ben}{\begin{eqnarray}}
\newcommand{\een}{\end{eqnarray}}
\newcommand{\beno}{\begin{eqnarray*}}
\newcommand{\eeno}{\end{eqnarray*}}
\newcommand{\me}{\mathrm{e}}

\newcommand{\sgn}{ {\rm sgn} }

\newtheorem{definition}{Definition}[section]

\newtheorem{thm}{Theorem}[section]
\newtheorem{lem}{Lemma}[section]
\newtheorem{rmk}{Remark}[section]

\newtheorem{prop}{Proposition}[section]
\renewcommand{\theequation}{\thesection.\arabic{equation}}

\def\sign{{\rm sign}}

\newdimen\eqjot \eqjot = 1\jot
\def\openupeq{\openup \the\eqjot}

\def\qaeq#1#2{{\def\\{&}\vcenter{\openupeq\halign{$\displaystyle
   ##\hfil$&&\hskip#1pt$\displaystyle##\hfil$\cr #2\cr}}}}

\def\qeq{\qaeq{20}}

\def\pofbox#1 #2$#3${\setbox0=\hbox{$#3$}\ht0=0pt\dp0=0pt\wd0=0pt\hskip-#1pt\raise#2pt\box0\hskip#1pt}

\begin{document}
\title[on the Cauchy problem for the modified CH equation with cubic nonlinearity]
{On the Cauchy problem for the integrable modified  Camassa-Holm  equation with cubic nonlinearity}
\author{Ying Fu}
\address{Ying Fu\newline
Department of Mathematics, Northwest University, Xi'an 710069,
P. R. China} \email{fuying@nwu.edu.cn}
\author{Guilong Gui}
\address{Guilong Gui\newline Department of Mathematics, Northwest University, Xi'an 710069,
P. R. China; \newline The Institute of Mathematical Sciences, The Chinese University of Hong Kong,  Hong Kong} \email{glgui@amss.ac.cn}
\author{Yue Liu}
\address{Yue Liu\newline
Department of Mathematics, University of Texas, Arlington, TX 76019; \newline  Department of Mathematics,
Ningbo University, P. R. China}
\email{yliu@uta.edu}
\author{Changzheng Qu}
\address{Changzheng Qu\newline
Department of Mathematics,
Ningbo University, Ningbo 315211,
P. R. China}  \email{quchangzheng@nbu.edu.cn}
\maketitle \numberwithin{equation}{section}
\begin{abstract} Considered in this paper is  the modified Camassa-Holm equation with cubic
nonlinearity, which is integrable  and admits the single peaked solitons and multi-peakon solutions. The  short-wave limit of this equation is known as  the
short-pulse equation. The main investigation is the Cauchy problem
of the modified Camassa-Holm equation with qualitative properties of
its solutions. It is firstly shown that the equation is locally
well-posed in a range of the Besov spaces. The blow-up scenario and
the lower bound of the maximal time of existence  are then
determined. A blow-up mechanism for solutions with certain initial
profiles is described in detail and nonexistence of the smooth
traveling wave solutions is also demonstrated.  In addition,  the
persistence properties of the strong solutions for the equation are
obtained.
\end{abstract}

\vskip 0.1cm


\noindent \small 2000 {\it Mathematics Subject Classification.} \
{\it Primary:} 35B30, 35G25.

\renewcommand{\theequation}{\thesection.\arabic{equation}}
\setcounter{equation}{0}

\section{Introduction}
In this paper, we are concerned with the following Cauchy  problem of the integrable  modified
Camassa-Holm  equation with cubic nonlinearity,
\begin{equation}\label{e1.1}
\left\{
 \begin{array}{ll}
\begin{split}
&m_t+(u^2-u^2_x)m_x+2u_xm^2+\gamma\,u_x=0, \quad  m=u-u_{xx}, \quad t>0,  \quad x \in \mathbb{R}, \\
&u(0,x)=u_0(x), \quad x \in \mathbb{R}.
 \end{split}
\end{array} \right.
\end{equation}
The equation in (\ref{e1.1}) was introduced by Fuchssteiner
\cite{fu} and Olver and Rosenau \cite{olv} (see also \cite{fok}) as
a new generalization of integrable system by implementing a simple
explicit algorithm based on the bi-Hamiltonian representation of the
classically integrable system. It also arises from a non-stretching
invariant curve flow in the two-dimensional Euclidean geometry
\cite{gloq}. In most cases, these new nonlinear systems are endowed
with nonlinear dispersion, and thus support non-smooth soliton-like
structures. It was  shown in \cite{qia} that the equation in
(\ref{e1.1}) admits the Lax-pair and the Cauchy problem (\ref{e1.1})
may be solved by the inverse scattering transform method. It was also found that the equation in (\ref{e1.1}) is related to the
short-pulse equation derived by Sch\"{a}fer and Wayne \cite{sw},
\begin{equation}\label{short}
v_{xt} = \frac{1}{3} (v^3)_{xx} + \gamma v,
\end{equation}
which is a model for the propagation of ultra-short light pulses in silica
optical fibers \cite{sw} and is also an approximation of nonlinear
wave packets in dispersive media in the limit of few cycles on the
ultra-short pulse scale \cite{cjsw}.

Actually, the short-pulse equation (\ref{short}) is a short-wave
limit of the equation in (\ref{e1.1}) by applying the following
scaling transformation \cite{gloq}
$$\qeq{x\longmapsto \epsilon\, x,\\t\longmapsto \epsilon^{-1} t,\\ u\longmapsto \epsilon^2 u}$$
where
$$u(t,x)= u_0(t,x) + \epsilon\, u_1(t,x)  + \epsilon^2\, u_2(t,x) + \>\cdots
$$
is expanded in powers of the  small parameter $\epsilon$. Then
$v=u_{0,x}(t,x)$ satisfies the short-pulse equation (\ref{short}).

The equation in (\ref{e1.1}) is formally  integrable and can be rewritten as
the bi-Hamiltonian form \cite{olv}, that is
\begin{equation*}
m_t=-((u^2-u^2_x)m)_x-\gamma u_x=J\frac{\delta H_0}{\delta
m}=K\frac{\delta H_1}{\delta m},
 \end{equation*}
where
\begin{equation*}
J=-\partial m\partial^{-1}m\partial-\frac{\gamma}{2}\partial,\quad
{\rm and} \qquad K=\partial^3-\partial,
\end{equation*} corresponding to  the Hamiltonian
\begin{equation*}
H_0=\int_{\mathbb{R}}mu \ dx,
\end{equation*}
and the Hamiltonian
\begin{equation*}
H_1=\frac1{4}\int_{\mathbb{R}} \left
(u^4+2u^2u^2_x-\frac1{3}u^4_x+2\gamma\,u^2 \right )\ dx.
\end{equation*}
It also admits  the Lax pair \cite{qia}, that is
\begin{align*}
\begin{pmatrix} \psi_1\\\psi_2\end{pmatrix}_x=U(m,\lambda)\begin{pmatrix}
\psi_1\\\psi_2\end{pmatrix}, \qquad \quad\begin{pmatrix}
\psi_1\\\psi_2\end{pmatrix}_t=V(m,u,\lambda)\begin{pmatrix}
\psi_1\\\psi_2\end{pmatrix},
\end{align*}
where
\begin{align*}
U(m,\lambda)=\frac12 \begin{pmatrix}-Q & \lambda
\,m\\-\lambda \,m &  Q\end{pmatrix}, \\
Q=Q(\lambda, \gamma)=\sqrt{1+\lambda^2\gamma}\,,
\end{align*}
and
\begin{align*}
V(m,u,&\lambda)=\\&{}-\frac{1}{2}\begin{pmatrix}\lambda^{-2}Q
+\frac1{2}Q(u^2-u^2_x) & -\lambda^{-1}(u-Q u_x)-\frac1{2}\lambda
(u^2-u^2_x)m \\  \lambda^{-1}(u+Q u_{x})+\frac{1}{2}\lambda
(u^2-u^2_x)m & -\lambda^{-2}Q-\frac{1}{2}Q(u^2-u^2_x)\end{pmatrix}.
\end{align*}

The Camassa-Holm (CH) equation \cite{cam, fuc} defined by
\begin{equation*}
m_t + u m_x + 2 u_x m + \gamma u_x  = 0,\;\; m=u-u_{xx}
\end{equation*}
has attracted much attention in the last twenty years because of its
interesting properties:  complete integrability, existence of peaked
solitons and multi-peakons \cite{cam, cao}, geometric formulations
\cite{cho1, coko, kou, mis} and the presence of breaking waves (i.e.
a solution that remains bounded while its slope becomes unbounded in
finite time) \cite{Co1,CE,CoEs2,CoEs3}. Note that the nonlinearity
in the CH equation is quadratic. In contrast to the integrable
modified KdV equation with a cubic nonlinearity, it is our interest
to find an integrable CH-type equations with a cubic nonlinearity.
Indeed, two integrable CH-type equations with cubic nonlinearity
have been discovered recently. One is the equation in (\ref{e1.1})
and the second one is the so-called Novikov equation \cite{nov}. The
integrability, peaked solitons, well-posedness and blow up phenomena
to the Novikov equation have been studied extensively, see the
references \cite{hon1, hon2, nov, tig1, tig2}, for example.

The goal of the present paper is to establish qualitative results
for the Cauchy problem (\ref{e1.1}).

We first study the local well-posedness  for the strong solutions to
the Cauchy problem (\ref{e1.1}) (see Theorem \ref{t3.1}). The proof
of the local well-posedness is inspired by the argument of
approximate solutions by Danchin \cite{Dan1} in the study of the
local well-posedness to the CH equation. However, one problematic
issue is that we here deal with a higher-order nonlinearity in the
Besov spaces, making the proof of several required nonlinear
estimates somewhat delicate. These difficulties are nevertheless
overcomed by carefully estimates for each iterative approximation of
solutions to the Cauchy problem (\ref{e1.1}).

With the local well-posedness
 obtained in hand, we then present a refined local well-posedness,
 i.e. local existence in the Besov space
$B^s_{2,1}$ with the critical index $s=\frac{5}{2}$ (see Theorem
\ref{t3.2}). Then a precise blow-up scenario (see Theorem
\ref{t4.2}) and a lower bound of the maximal time of existence (see
Theorem \ref{t4.3}) are obtained.

Blow-up in finite time depends on strong nonlinear dispersion
usually and makes, of course, the analysis more challenging in our
case with higher nonlinearities. It is known that  a solution of the
Camassa-Holm equation, which can be considered as the transport
equation, blows up in finite time when its slope $ u_x $ is
unbounded from below. This idea is expected to be applied to the
modified CH equation in \eqref{e1.1}, since it can be written as  a
transport equation in terms of $m$ along the flow generated by
$u^2-u_x^2$, that is
\begin{equation} \label{tmch}
m_t+(u^2-u^2_x)m_x=-2u_x m^2-\gamma u_x.
\end{equation}
Generally speaking, the transport equation theory ensures that, if the
slope
\begin{equation}\label{slope}
 (u^2 - u^2_x)_x = 2 u_x m
\end{equation}
is bounded, the solution will remain regular and, therefore, can not
blow up in finite time.  In view of this property, together with the
Sobolev embedding theorem, it can be shown that the solution blows
up in finite time if and only if the slope in \eqref{slope} is
unbounded from below. Thus to prevent the solution from blow-up in
finite time, the main issue is that it is impossible to control the
bound of $ u_x m $ in \eqref{slope} in terms of the $ H^1$-norm of
the solution unless a higher, positive conserved quantity involved
in $ H^3 $-norm of the solution can be found. To overcome this
difficulty, we may regard the evolution equation \eqref{tmch} in
terms of the quantity \eqref{slope} being transported along the flow
generated by $ u^2 - u^2_x $. Then blow-up result can be established
by using the global conservative property of the potential density
$m$ along the characteristics. This new idea was used in \cite
{gloq} in the case $ \gamma = 0. $ Inspired by this method, we are
able to improve the blow up result in \cite{gloq} by using the
conservation quantity $ \displaystyle I_0 = \int_{\mathbb {R}} u(0,
x) dx = \int_{\mathbb{R}} u(t, x) dx $
 (see Theorem \ref{t5-1}). Moreover, with some initial data decaying exponentially, we prove
that the solution to the initial value problem (\ref{e1.1}) is also
decaying exponentially, i.e. the persistence property.

As mentioned above, it  is well known  that the CH equation has the
peakons \cite{cam}, which are shown to be orbitally stable in the
intriguing papers \cite{CoM, CoS1}. Stability of the periodic
peakons of the CH equation can be found in \cite{len}. So it is of
interest to identify traveling-wave solutions of the equation in
(\ref{e1.1}). Indeed, Gui-Liu-Olver-Qu in \cite{gloq} found that the
equation in (\ref{e1.1}) with $ \gamma = 0 $  has single peakons
given by
\begin{align*}
u_c(t, x)= \sqrt{\frac{3c}{2}}\,e^{-|x-ct|},  \quad c > 0
\end{align*}
and multi-peakons. In particular, the two-peakons can be given
explicitly by
\begin{equation*}
\begin{split}
u(t, x) =  & \sqrt{\frac 32 c_1} \,\exp\left\{- \left |x-c_1 t -\frac
{3\sqrt{c_1c_2}}{c_1-c_2}e^{(c_1-c_2)t} \right |\right \}\\&
+\sqrt{\frac 32 c_2} \, \exp\left\{- \left |x-c_2 t-\frac
{3\sqrt{c_1c_2}}{c_1-c_2}e^{(c_1-c_2)t} \right |\right \}, \quad 0 < c_1 < c_2.
\end{split}
\end{equation*}
As a part of the present paper, we are able to show that the equation in (\ref{e1.1}) with $ \gamma = 0 $  does
not have any nontrivial smooth traveling-wave solutions.

The rest of  the paper is organized as follows. In Section 2, some
preliminary properties, which will be used later, are presented. The
local well-posedness in the Besov spaces is established in Section
3. In Section 4, a blow-up scenario and a lower bound of the maximal
existence time of (\ref{e1.1}) will be derived. A new blow-up
mechanism is described and some blow-up data are determined  in
Section 5. Non-existence of smooth traveling waves for $ \gamma = 0
$ is demonstrated in Section 6.   Section 7 is devoted to  the
persistence properties of the initial-value problem (\ref{e1.1}).

\vskip 0.1cm

\noindent {\it Notation}. In the following, for a given Banach space
$Z$, we denote its norm by $\|\cdot\|_Z$. Since all space of
functions are over $\mathbb{R}$, for simplicity, we drop
$\mathbb{R}$ in our notations of function spaces if there is no
ambiguity. We denote ${\cF}u$ or $\hat{u}$ the Fourier transform of
the function $u$.

\renewcommand{\theequation}{\thesection.\arabic{equation}}
\setcounter{equation}{0}

\section{Preliminaries}

For the convenience of the reader, we shall recall some basic facts on the Littlewood-Paley
theory and the transport equations theory, one may check \cite{BCD, Che, Dan1, Dan2} for more details.

\begin{prop}\label{p2.1}\cite{BCD, Che} (Littlewood-Paley decomposition) Let
${\cB}\eqdefa\{\xi\in\mathbb{R},\ |\xi|\leq \frac{4}{3}\}$  and
${\cC}\eqdefa\{\xi\in\mathbb{R},\ \frac{3}{4}\leq |\xi|\leq
\frac{8}{3}\}.$ There exist two radial functions $\chi\in
C_c^\infty({\cB})$ and $ \varphi\in C_c^\infty({\cC})$ such that
\begin{equation*}
\chi(\xi)+ \sum_{q  \geq 0}\varphi(2^{-q}\xi)=1,\quad \forall \ \ \
 \xi\in \mathbb{R}^d ,
\end{equation*}
\begin{equation*}
|q-q^{\prime}|\geq 2 \Rightarrow \Supp \varphi(2^{-q}\cdot)\cap
\Supp
 \varphi(2^{-q^{\prime}}\cdot) = \varnothing,
\end{equation*}
\begin{equation*}
q\geq 1 \Rightarrow \Supp \chi(\cdot)\cap \Supp
 \varphi(2^{-q}\cdot) = \varnothing,
\end{equation*}
and
\begin{equation*}
\frac{1}{3}\leq \chi(\xi)^2+\sum_{q \geq 0}\varphi((2^{-q}\xi))^2
\leq 1, \quad \forall \ \ \
 \xi\in \mathbb{R}^d .
\end{equation*}
\end{prop}

Furthermore, let $h \eqdefa {\cF}^{-1}\varphi$ and $\tilde{h}
\eqdefa {\cF}^{-1}\chi.$ Then the dyadic operators $\Delta_q$ and
$S_q$ can be defined as follows
\begin{equation*}\begin{split}
&\Delta_qf \eqdefa
\varphi(2^{-q}D)f=2^{qd}\int_{\R^d}h(2^qy)f(x-y)dy
\quad\mbox{for}\quad
q\geq 0, \nonumber\\
&S_qf \eqdefa \chi(2^{-q}D)f=\sum_{-1\leq k\le
q-1}\Delta_kf=2^{qd}\int_{\mathbb{R}^d}\tilde{h}(2^qy)f(x-y)dy, \nonumber\\
&\D_{-1}f \eqdefa S_{0}f \, \, \mbox{and} \, \, \D_{q}f \eqdefa 0
\quad \,\mbox{for} \quad\, q \leq -2.
\end{split}
\end{equation*}

\begin{definition}\label{d2.1}\cite{BCD, Che} (Besov space)
Let $s\in \mathbb{R}, 1\le p,r\le\infty.$ The  inhomogenous Besov
space $B^s_{p,r}(\mathbb{R}^{d})$ ($B^s_{p,r}$ for short) is defined
by
\begin{equation*}
B^s_{p,r}\eqdefa\{f\in {\mathcal S}^{\prime}(\mathbb{R}^{d}); \quad
\|f\|_{B^s_{p,r}}<\infty\},
\end{equation*}
where
\begin{equation*}
\|f\|_{B^s_{p,r}}\eqdefa\left\{\begin{array}{l}
\displaystyle\bigg(\sum_{q \in \mathbb{Z}} 2^{q s r}\|\Delta_q
f\|_{L^p}^r\bigg)^{\frac 1
r},\quad \hbox{for}\quad r<\infty,\\
\displaystyle\sup_{q \in \mathbb{Z}}2^{q s}\|\Delta_q f\|_{L^p},
\quad \quad \ \quad \hbox{ for} \quad r=\infty.
\end{array}\right.
\end{equation*}
If $s=\infty,$ $B^{\infty}_{p, r} \eqdefa \bigcap_{s \in \mathbb{R}}
B^{s}_{p, r} .$
\end{definition}

\begin{prop}\label{p2.2}\cite{BCD, Dan1,Dan2}
The following properties hold.

i) Density: if $p, \, r < \infty,$ then $\mathcal{S}(\mathbb{R}^d)$
is dense in $B^s_{p, r}(\mathbb{R}^d).$

ii) Sobolev embeddings: if $p_1 \leq p_2$ and   $r_1 \leq r_2,$ then
$B^s_{p_1, r_1} \hookrightarrow
B^{s-d(\frac{1}{p_1}-\frac{1}{p_2})}_{p_2, r_2}.$ If $s_1 <s_2,$ $1
\leq p \leq +\infty$ and $1 \leq r_1, \, r_2 \leq +\infty,$ then the
embedding $B^{s_2}_{p, r_2} \hookrightarrow B^{s_1}_{p, r_1}$ is
locally compact.

iii) Algebraic properties: for $s>0,$ $B^s_{p, r}\cap L^{\infty}$ is
an algebra. Moreover, ($B^s_{p, r}$ is an algebra)
$\Longleftrightarrow$ ($B^s_{p, r} \hookrightarrow L^{\infty}$)
$\Longleftrightarrow$ ($s > \frac{d}{p}$ or ($s \geq \frac{d}{p}$
and $r=1$)).

iv) Fatou property: if $(u^{(n)})_{n\in \mathbb{N}}$ is a bounded
sequence of $B^s_{p, r}$ which tends to $u$ in
${\mathcal{S}}^{\prime},$  then $u \in B^s_{ p,r}$ and
\begin{equation*}
\|u\|_{B^s_{ p,r}} \leq \lim\inf_{n \rightarrow \infty}
\|u^{(n)}\|_{B^s_{ p,r}}.
\end{equation*}

v) Complex interpolation: if $u \in B^s_{ p,r} \cap B^{\tilde{s}}_{
p,r}$ and $\theta \in [0, 1], \, 1 \leq p, r \leq \infty,$ then $u
\in  B^{\theta s+ (1-\theta)\tilde{s}}_{ p,r}$ and $\|u\|_{B^{\theta
s+ (1-\theta)\tilde{s}}_{ p,r}} \leq \|u\|_{B^s_{ p,r}}^{\theta}
\|u\|_{B^{\tilde{s}}_{ p,r}}^{1-\theta}.$

vi) Let $m \in \mathbb{R}$ and $f$ be a $S^m$-multiplier (that is,
$f: \mathbb{R}^{d} \rightarrow \mathbb{R}$ is smooth and satisfies
that for all multi-index $\alpha,$ there exists a constant
$C_{\alpha}$ such that for any $\xi \in \mathbb{R}^{d}$,
$|\partial^{\alpha} f(\xi)| \leq C_{\alpha}
(1+|\xi|)^{m-|\alpha|}.$) Then for all $s \in \mathbb{R}$ and $1
\leq p, r \leq \infty,$ the operator $f(D)$ is continuous from $
B^s_{ p,r} $ to $ B^{s-m}_{ p,r} .$
\end{prop}

\begin{lem}\label{l2.2}\cite{BCD, Dan1,Dan2}
Suppose that $(p, r)\in [1, +\infty]^2$ and $s >-\frac{d}{p}.$ Let
$v$ be a vector field such that $\nabla v$ belongs to $L^1([0, T];
B^{s-1}_{p, r})$ if $s >1+\frac{d}{p}$ or to  $L^1([0, T];
B^{\frac{d}{p}}_{p, r}\cap L^{\infty})$ otherwise. Suppose also that
$f_{0}\in B^{s}_{p, r}, \, F \in L^1([0, T]; B^{s}_{p, r})$ and that
$f \in L^{\infty}([0, T]; B^{s}_{p, r})\cap C([0, T];
\mathcal{S}^{\prime})$ solves the d-dimensional linear transport
equations
\begin{equation*}(T)~~~~~~~~
\begin{cases}
\partial_{t}f + v \cdot \nabla f= F,\\
f|_{t=0}=f_0.
\end{cases}
\end{equation*}
Then there exists a constant $C$ depending only on $s, \, p$ and $d$
such that the following statements hold:

1) If $r=1$ or $s \neq 1+\frac{d}{p},$ then
\begin{equation*}
\|f\|_{B^{s}_{p, r}}\leq \|f_0\|_{B^{s}_{p, r}}+\int_{0}^{t}
\|F(\tau)\|_{B^{s}_{p, r}}d\tau+C\int_{0}^{t}
V'(\tau)\|f(\tau)\|_{B^{s}_{p, r}}d\tau,
\end{equation*}
or
\begin{equation}\label{2.2}
\|f\|_{B^{s}_{p, r}}\leq e^{CV(t)}\left( \|f_0\|_{B^{s}_{p,
r}}+\int_{0}^{t} e^{-CV(\tau)}\|F(\tau)\|_{B^{s}_{p, r}}d\tau
 \right)
\end{equation}
hold, where  $V(t)= \int_{0}^{t} \|\nabla
v(\tau)\|_{B^{\frac{d}{p}}_{p, r} \cap L^{\infty}}d\tau$ if $s < 1+
\frac{d}{p}$ and $V(t)= \int_{0}^{t} \|\nabla v(\tau)\|_{B^{s-1}_{p,
r} }d\tau$ else.

2) If  $s \leq 1+\frac{d}{p}$ and, in addition, $\nabla f_0 \in
L^{\infty},$ $\nabla f \in L^{\infty}([0, T] \times \mathbb{R}^{d})$
and $\nabla F \in L^{1}([0, T]; L^{\infty}),$ then
\begin{equation*}
\begin{split}
 \|f(t)&\|_{B^{s}_{p, r}}+\|\nabla f(t)\|_{L^{\infty}}\\
  \leq & e^{CV(t)}\left(
\|f_0\|_{B^{s}_{p, r}}+\|\nabla f_0\|_{L^{\infty}} +\int_{0}^{t}
e^{-CV(\tau)}(\|F(\tau)\|_{B^{s}_{p, r}}+\|\nabla
F(\tau)\|_{L^{\infty}})d\tau
 \right)
\end{split}
\end{equation*}
with  $V(t)= \int_{0}^{t} \|\nabla v(\tau)\|_{B^{\frac{d}{p}}_{p, r}
\cap L^{\infty}}d\tau.$

3) If $f=v,$ then for all $s >0,$ the estimate (\ref{2.2}) holds
with $V(t)=\int_{0}^{t}\|\partial_x u(\tau)\|_{L^{\infty}}d\tau .$

4) If $r <+\infty,$ then $f \in C([0, T];B^{s}_{p, r} ).$ If $r
=+\infty,$ then $f \in C([0, T]; B^{s^{\prime}}_{p, 1} )$ for all
$s^{\prime} <s.$
\end{lem}

\begin{lem}\label{l2.3}\cite{Dan2}
Let $(p, p_1,  r)\in [1, +\infty]^3.$ Assume that  $s
>-d\min\{\frac{1}{p_1}, \frac{1}{p^{\prime}}\}$ with $p^{\prime} \eqdefa (1-\frac{1}{p})^{-1}.$
 Let $f_0\in B^{s}_{p, r}$ and $F \in L^1([0, T]; B^{s}_{p, r}).$
Let $v$ be a time dependent vector field such  that $v \in
L^{\rho}([0, T]; B^{-M}_{\infty, \infty})$ for some $\rho >1,$ $M>0$
and $\nabla v \in L^{1}([0, T]; B^{\frac{d}{p_1}}_{p_1, \infty}\cap
L^{\infty} )$ if $s < 1+\frac{d}{p_1},$ and $\nabla v \in L^{1}([0,
T]; B^{s-1}_{p_1, r} )$ if $s > 1+\frac{d}{p_1}$ or $s =
1+\frac{d}{p_1}$ and $r=1.$ Then the transport equations (T) has a
unique solution $f \in L^{\infty}([0, T]; B^s_{p, r}) \cap
(\cap_{s^{\prime}<s}C([0, T]; B^{s^{\prime}}_{p, 1}))$ and the
inequalities in Lemma \ref{l2.2} hold true. If, moreover, $r <
\infty,$ then we have $f \in C([0, T];B^s_{p, r}).$
\end{lem}

\begin{lem}\label{l2.4}\cite{Che} (1-D Moser-type estimates) Assume
that $1 \leq p, \, r \leq +\infty,$ the following estimates hold:

(i)$\;\;$ for $s>0,$  $\|fg\|_{B^{s}_{p, r}}\leq C (\|f\|_{B^{s}_{p,
r}}\|g\|_{L^{\infty}}+ \|g\|_{B^{s}_{p, r}}\|f\|_{L^{\infty}});$

(ii)$\;\;$ for $s_1 \leq \frac{1}{p},\, s_2>\frac{1}{p}$ ($ s_2\geq
\frac{1}{p}$ if $r=1$) and $s_1+s_2>0,$
\begin{equation*}
\|fg\|_{B^{s_1}_{p, r}}\leq C \|f\|_{B^{s_1}_{p,
r}}\|g\|_{B^{s_2}_{p, r}},
\end{equation*}
where  the constant $C$ is independent of $f$ and $g$.
\end{lem}




\begin{lem}\label{l2.5}\cite{Dan2}
Denote $\bar{\mathbb{N}}=\mathbb{N}\cup{\infty}$. Let
$(v^{(n)})_{n\in\bar{\mathbb{N}}}$ be a sequence of functions
belonging to $C([0,T];B^{\frac1{2}}_{2,1})$. Assume that $v^{(n)}$
is the solution to
\begin{equation}\label{e2.2}
\begin{cases}
\partial_{t}v^{(n)} + a^{(n)}\partial_{x}v^{(n)}= f,\\
v^{(n)}|_{t=0}=v_0
\end{cases}
\end{equation}
with $v_0\in B^{\frac1{2}}_{2,1},f\in L^1(0,T;B^{\frac1{2}}_{2,1})$
and that, for some $\beta\in L^1(0,T)$,
$$\sup\limits_{n\in \mathbb{N}}\|\partial_x a^{(n)}(t)\|_{B^{\frac1{2}}_{2,1}}\leq\beta(t).$$
If in addition $a^{(n)}$ tends to $a^{\infty}$ in
$L^1(0,T;B^{\frac1{2}}_{2,1})$ then $v^{(n)}$ tends to $v^{\infty}$
in $C(0,T;B^{\frac1{2}}_{2,1})$.
\end{lem}

Next we reformulate the Cauchy problem (\ref{e1.1}) in a more
convenient form. Note that the equation in (\ref{e1.1}) is
equivalent to the following one.
\begin{equation*}
u_t-u_{xxt}+3u^2u_x-4uu_xu_{xx}+u^2_xu_{xxx}+2u_xu^2_{xx}-u^2u_{xxx}-u^3_x+\gamma\,u_x=0.
\end{equation*}
Applying the operator $(1-\partial_x^2)^{-1}$ to both sides of the
above equation, we obtain
\begin{equation}\label{eqns-2}
u_t+ \left (u^2-\dfrac1{3}u^2_x \right
)u_x+\partial_x(1-\partial_x^2)^{-1} \left (\dfrac{2}{3}u^3+uu^2_x
\right )
+(1-\partial_x^2)^{-1}\left(\dfrac{u^3_x}{3}+\gamma\,u_x\right)=0,
\end{equation}
which enables us to define the weak solution of the Cauchy problem \eqref{e1.1}.
\renewcommand{\theequation}{\thesection.\arabic{equation}}
\setcounter{equation}{0}
\section{Local well-posedness}
\subsection{Local existence}
In this section, we shall discuss the local well-posedness of the
Cauchy  problem (\ref{e1.1}). At first, we introduce the following
spaces.

\begin{definition}\label{d3.1}
For $T > 0, \, s \in \mathbb{R}$ and $1 \leq p \leq
+\infty,$ we set
\begin{equation*}
\begin{split}
&E^s_{p, r}(T) \eqdefa C([0, T]; B^s_{p, r})\cap C^{1}([0, T];
B^{s-1}_{p, r}) \, \quad \mbox{if} \quad
\, r< +\infty,\\
&E^s_{p, \infty}(T) \eqdefa L^{\infty}([0, T]; B^s_{p, \infty})\cap
Lip([0, T]; B^{s-1}_{p, \infty})
\end{split}
\end{equation*}
and $E^s_{p, r} \eqdefa \cap_{T >0}E^s_{p, r}(T).$
\end{definition}

The result of the local well-posedness in the Besov space may now be
enunciated.

\begin{thm}\label{t3.1}
Suppose that $1 \leq p, \, r \leq +\infty $, $s>\max
\{2+\frac{1}{p}, \frac{5}{2}\}$ and $u_0 \in B^{s}_{p, r}$. Then
there exists a time $T> 0$ such that the initial-value problem
(\ref{e1.1})
 has a unique solution $u \in E^s_{p, r}(T)
,$ and the map $u_0 \mapsto u$ is continuous from a neighborhood of
$u_0$ in $B^s_{p, r} $ into
\begin{equation*}
C([0, T]; B^{s^{\prime}}_{p, r}) \cap C^1([0, T];
B^{s^{\prime}-1}_{p, r})
\end{equation*} for every
$s^{\prime} < s$ when $r =+\infty$ and  $s^{\prime} =s$ whereas $r
<+\infty.$
\end{thm}
\begin{rmk}
When $p=r=2,$ the Besov space $B^{s}_{p, r}$ coincides with the
Sobolev space $H^{s}.$ Theorem \ref{t3.1} implies that under the
condition $u_0 \in H^{s} $ with $s> 5/2,$ we can obtain the local
well-posedness  for the initial-value problem (\ref{e1.1}).
\end{rmk}
\begin{rmk} As in Remark 4.1 in \cite{gloq}, the existence time for the initial-value problem (\ref{e1.1}) may be
chosen independently of $s$ in the following sense. If
\begin{equation*}
u \in C([0, T]; H^s) \cap C^1([0, T]; H^{s-1})
\end{equation*} is the solution of the initial-value problem
(\ref{e1.1}) with initial data $u_0 \in H^r$ for some $r
>5/2, \, r \neq s,$ then
\begin{equation*}
u\in C([0, T]; H^r) \cap C^1([0, T]; H^{r-1})
\end{equation*}
 with the same time $T$. In particular, if $u_0 \in H^{\infty},$ then $u \in C([0, T];
H^{\infty})$.
\end{rmk}
\begin{rmk}\label{rmk-conservation-1}
For a strong solution $m = u - u_{xx}$ in Theorem \ref{t3.1}, if, in
addition, the initial data $u_0 \in L^1$,
 then the following three functionals are conserved:
\begin{equation}\label{cons-0-1}
I_0=\int_{\mathbb{R}} u(t) \, dx, \quad I_1=\int_{\mathbb{R}}(u^2+u_x^2)\, dx,\quad I_2=\int_{\mathbb{R}}
\big(u^4+2u^2u^2_x-\frac1{3}u^4_x+2\gamma\,u^2 \big)\ dx.
\end{equation}
\end{rmk}

Under the assumptions in Theorem \ref{t3.1} (especially $p=r=2$), we introduce the flow generated by
$u^2-u^2_x$:
\begin{equation}\label{flow-1}
\left\{
\begin{aligned}
\frac{d\,q(t, x)}{d\,t}&=(u^2-u^2_x)(t,q(t,x)), \\
q(0,x)&=x,
 \end{aligned}
\right.  \qquad x\in \mathbb{R},\quad t\in[0,T),
\end{equation}

If $\gamma=0$, then it is easy to check that \cite{gloq}
\begin{equation*}\label{positive-1}
m(t,q(t,x))q_{x}(t,x) =m_{0}(x),\quad \hbox{for all} \quad  (t,x)\in
[0,T)\times\mathbb{R}.
\end{equation*}

\begin{rmk}\label{rmk-4-1}
Note that from the above flow, it follows that \cite{gloq}
\begin{equation*}
q_{x}(t,x)= \left ( 2\int_{0}^{t}(mu_{x})(s,q(s,x))ds \right ) >
0,\, \quad \hbox{for all} \quad  (t,x)\in [0,T)\times \mathbb{R}.
\end{equation*}
In view of the above
conservation law, we deduce that: if $u_0(x)$ has compact support in
$x$ in the interval $[a,b]$, then so does $m(t,\cdot)$ in the
corresponding interval $[q(t,a),q(t,b)]$. Moreover, if $m_0 =(1-\partial_x^2) u_0$
does not change sign, then $m(t,x)$ will not change sign for any $ t \in [0, T)$. On the other hand, the
$L^{\infty}$-norm of any function $v(t,\cdot)\in L^{\infty}$ is preserved under the family of diffeomorphisms
$q(t,\cdot)$, that is,
\begin{equation*}
\|v(t,\cdot)\|_{L^{\infty}}=\|v(t,q(t,\cdot))\|_{L^{\infty}},\quad
t\in[0,T).
\end{equation*}
\end{rmk}

In the following, we denote $C>0$ a generic constant only depending
on $p,\; r,\; s$. Uniqueness and continuity with respect to the
initial data are an immediate consequence of the following result.

\begin{prop}\label{p3.1}
 Let $1 \leq p, \, r \leq +\infty $ and $s > \max\{2+\frac{1}{p}, \, \frac{5}{2}\}.$
 Let $u^{(1)}, \,u^{(2)}$ be two given solutions of the initial-value problem
(\ref{e1.1}) with the initial data $u^{(1)}_{0}, \, u^{(2)}_{0} \in
B^{s}_{p, r}$ satisfying $u^{(1)}, \, u^{(2)} \in L^{\infty}([0, T];
B^{s}_{p, r})\cap C([0, T]; \mathcal{S}^{\prime})$. Then for every
$t \in [0, T]:$
\begin{equation}\label{3.1}
\begin{split}
&\|(u^{(1)}-u^{(2)})(t)\|_{B^{s-1}_{p, r}}\\
&\leq \|u^{(1)}_0-u^{(2)}_0\|_{B^{s-1}_{p, r}} \exp \left
\{C\int_{0}^{t}\left(\|u^{(1)}(\tau)\|^2_{B^{s}_{p,
r}}+\|u^{(2)}(\tau)\|^2_{B^{s}_{p, r}}+|\gamma|\right)\;d\tau
\right\}.
\end{split}
\end{equation}
\end{prop}

\begin{proof} Denote $u^{(12)} \eqdefa u^{(2)}-u^{(1)}$. It is obvious that
$$u^{(12)} \in
L^{\infty}([0, T]; B^{s}_{p, r})\cap C([0, T];
\mathcal{S}^{\prime}),$$
which along with the equivalent formulation \eqref{eqns-2} of \eqref{e1.1} implies that $u^{(12)} \in C([0, T];
B^{s-1}_{p, r})$ and $u^{(12)}$ solves the transport equation
\begin{equation}\label{e-d-1}
\begin{split}
\partial_t u^{(12)}+\left[(u^{(1)})^2-\frac{1}{3}((u^{(1)}_x)^2+u^{(1)}_xu^{(2)}_x+(u^{(2)}_x)^2)\right]\pa_x u^{(12)}=f(u^{(12)}, u^{(1)}, u^{(2)})
\end{split}
\end{equation}
with
\begin{equation}\label{r-d-1}
\begin{split}
f(u^{(12)}, u^{(1)}, u^{(2)})=&-(1-\partial_x^2)^{-1}\bigg(\frac{1}{3}((u^{(1)}_x)^2+u^{(1)}_xu^{(2)}_x
+(u^{(2)}_x)^2)\, u^{(12)}_x+ \gamma \, u^{(12)}_x\bigg)\\
&-(u^{(1)}+u^{(2)})u^{(2)}_x \,
u^{(12)}-\partial_x(1-\partial_x^2)^{-1}\bigg(\frac{2}{3}((u^{(1)})^2+u^{(1)}u^{(2)}\\
&+(u^{(2)})^2)\, u^{(12)}+
u^{(1)}(u^{(1)}_x+u^{(2)}_x) \,u^{(12)}_x +(u^{(2)}_x)^2\, u^{(12)}\bigg).
\end{split}
\end{equation}
Thanks to the transport theory
in Lemma \ref{l2.2}, one gets
\begin{equation}\label{e-d-2}
\begin{split}
\|u^{(12)}(t)\|_{B^{s-1}_{p,r}}\leq & C\int_0 ^t\left
\|(u^{(1)})^2-\frac{1}{3}((u^{(1)}_x)^2+u^{(1)}_xu^{(2)}_x+(u^{(2)}_x)^2)\right
\|_{B^{s-1}_{p,r}}\,\|u^{(12)}(\tau)\|_{B^{s-1}_{p,r}}\,
d\tau\\
&+ \int_0^t\|f(u^{(12)}, u^{(1)}, u^{(2)})(\tau)\|_{B^{s-1}_{p,r}}\,d\tau+\|u^{(12)}(0)\|_{B^{s-1}_{p,r}}.
\end{split}
\end{equation}
Applying the product law in the Besov spaces, we have
\begin{equation*}
\begin{split}
\left\|(u^{(1)})^2-\frac{1}{3}((u^{(1)}_x)^2+u^{(1)}_xu^{(2)}_x+(u^{(2)}_x)^2)
\right\|_{B^{s-1}_{p,r}} \leq C (\|u^{(1)}\|_{B^{s}_{p,r}}
^2+\|u^{(2)}\|_{B^{s}_{p,r}}^2).
\end{split}
\end{equation*}
Similarly, one gets
\begin{equation*}
\begin{split}
&\left\|(1-\partial_x^2)^{-1} \bigg(\frac{1}{3}((u^{(1)}_x)^2+u^{(1)}_xu^{(2)}_x+(u^{(2)}_x)^2)\, u^{(12)}_x
+ \gamma \, u^{(12)}_x\bigg)\right\|_{B^{s-1}_{p,r}}\\
&\leq C\left\|\frac{1}{3}((u^{(1)}_x)^2+u^{(1)}_xu^{(2)}_x+(u^{(2)}_x)^2)\, u^{(12)}_x+ \gamma \, u^{(12)}_x\right\|_{B^{s-2}_{p,r}} \\
&\leq
C\bigg(\|u^{(1)}\|_{B^{s-1}_{p,r}}^2+\|u^{(1)}\|_{B^{s-1}_{p,r}} \|u^{(2)}\|_{B^{s-1}_{p,r}}
+\|u^{(2)}\|_{B^{s-1}_{p,r}}^2+|\gamma|\bigg)\,
\|u^{(12)}\|_{B^{s-1}_{p,r}},
\end{split}
\end{equation*}
$$
\|(u^{(1)}+u^{(2)})u^{(2)}_x \, u^{(12)}\|_{B^{s-1}_{p,r}}  \leq C(\|u^{(1)}\|_{B^{s-1}_{p,r}} +\|u^{(2)}\|_{B^{s-1}_{p,r}} )
\|u^{(2)}\|_{B^{s}_{p,r}} \| u^{(12)}\|_{B^{s-1}_{p,r}} ,
$$
and
\begin{equation*}
\begin{split}
&\big\|\partial_x (1-\partial_x^2)^{-1}  \bigg(\frac{2}{3}((u^{(1)})^2+u^{(1)}u^{(2)}+(u^{(2)})^2)\, u^{(12)}\\
&\qquad\qquad \qquad +u^{(1)}(u^{(1)}_x+u^{(2)}_x) \,u^{(12)}_x +(u^{(2)}_x)^2\, u^{(12)}\bigg)\big\|_{B^{s-1}_{p,r}}\\
&\leq C\left\|\frac{2}{3}((u^{(1)})^2+u^{(1)}u^{(2)}+(u^{(2)})^2)\, u^{(12)}+u^{(1)}(u^{(1)}_x+u^{(2)}_x) \,u^{(12)}_x
 +(u^{(2)}_x)^2\, u^{(12)}\right\|_{B^{s-2}_{p,r}}\\
&\leq
C\bigg(\|u^{(1)}\|_{B^{s-1}_{p,r}}^2+\|u^{(2)}\|_{B^{s-1}_{p,r}}^2\bigg)\,
 \|u^{(12)}\|_{B^{s-1}_{p,r}},
\end{split}
\end{equation*}
which leads to
\begin{equation*}
\begin{split}
&\|f(u^{(12)}, u^{(1)}, u^{(2)})\|_{B^{s-1}_{p,r}}\leq
C\bigg(\|u^{(1)}\|_{B^{s}_{p,r}}^2+\|u^{(2)}\|_{B^{s}_{p,r}}^2+|\gamma|\bigg)\,
\|w\|_{B^{s-1}_{p,r}}.
\end{split}
\end{equation*}
Hence, one obtains from \eqref{e-d-2} that
\begin{equation*}
\begin{split}
&\|u^{(12)}(t)\|_{B^{s-1}_{p,r}}\\
& \leq
\|u^{(12)}(0)\|_{B^{s-1}_{p,r}}
+C\int_{0}^{t}\bigg(\|u^{(1)}(\tau)\|_{B^{s}_{p,r}}^2+\|u^{(2)}(\tau)\|_{B^{s}_{p,r}}^2
+|\gamma|\bigg) \|u^{(12)}(\tau)\|_{B^{s-1}_{p,r}}\, d\tau,
\end{split}
\end{equation*}
and then applying Gronwall's inequality, we reach \eqref{3.1}.
\end{proof}
Now let us start the proof of Theorem \ref{t3.1}, which is motivated
by the proof of local existence theorem about the Camassa-Holm
equation in \cite{Dan1}. Firstly, we shall use the classical
Friedrichs regularization method to construct the approximate
solutions to the Cauchy problem  (\ref{e1.1}).

\begin{lem}\label{l3.1}
Let $u_0,\, p, \, r $ and $ s $ be as in the statement of Theorem
\ref{t3.1}. Assume that $u^{(0)}  := 0.$ There exists a sequence of
smooth functions $(u^{(n)})_{n \in \mathbb{N}} \in C(\mathbb{R}^{+};
B^{\infty}_{p, r})$ solving the following linear transport equation
by induction:
\begin{equation}\label{3.3}(T_n)~~~~~~
\begin{cases}
\left\{\pa_t +\left[(u^{(n)})^2-(u^{(n)}_x)^2\right]
\pa_x\right\}m^{(n+1)}=-2u^{(n)}_x(m^{(n)})^2-\gamma u^{(n)}_x,\quad
& t> 0,\;x\in \mathbb{R},\\
u^{(n+1)}|_{t=0} = u_{0}^{(n+1)}(x)=S_{n+1}u_0,& x\in \mathbb{R}.
\end{cases}
\end{equation}
Moreover, there exists a $T>0$ such that the solutions satisfying
the following properties:

(i) $\;\;(u^{(n)})_{n \in \mathbb{N}}$ is uniformly bounded in $E^{s}_{p,
r}(T)$.

(ii) $\;\;(u^{(n)})_{n \in \mathbb{N}}$ is a Cauchy sequence in
$C([0, T]; B^{s-1}_{p, r})$.
\end{lem}

\begin{proof}
Since all data $S_{n+1}u_0$ belongs to $B^{\infty}_{p, r},$ Lemma
\ref{l2.3} enables us to show by induction that for all $n \in
\mathbb{N},$ the  equation $(T_n)$ has a global solution which
belongs to $C(\mathbb{R}; B^{\infty}_{p, r}).$ Applying Lemma
\ref{l2.2} to $(T_n)$, we get for all $n \in \mathbb{ N} :$
\begin{equation}\label{3.4-a}
\begin{split}
&\|m^{(n+1)}(t)\|_{B^{s-2}_{p, r}} \leq
e^{C\int_0^t\left\|\left[(u^{(n)})^2-(u^{(n)}_x)^2\right](\tau)\right\|_{B^{s-1}_{p,
r}}d\tau}\|S_{n+1}u_0\|_{B^{s}_{p,
r}}\\&+C\int_0^t
e^{C\int^{t}_{\tau}\left\|\left[(u^{(n)})^2-(u^{(n)}_x)^2\right](\tau')\right\|_{B^{s-1}_{p,
r}}d\tau^{\prime}} \|2u^{(n)}_x(m^{(n)})^2+\gamma u^{(n)}_x(\tau)\|_{B^{s-2}_{p,
r}}\;d\tau.
\end{split}
\end{equation}
Thanks to the product law in Besov spaces, one has
\begin{equation*}
\begin{split}
&\|(u^{(n)})^2-(u^{(n)}_x)^2\|_{B^{s-1}_{p,
r}} \leq C \|u^{(n)}\|_{B^{s}_{p,
r}}^2, \\
&\|2u^{(n)}_x(m^{(n)})^2+\gamma u^{(n)}_x\|_{B^{s-2}_{p,
r}}\leq C \left(\|u^{(n)}\|^3_{B^{s}_{p,
r}}+\|u^{(n)}\|_{B^{s}_{p, r}}\right),
\end{split}
\end{equation*}
which along with \eqref{3.4-a} leads to
\label{3.4}
\begin{equation}\label{3.4}
\begin{split}
\|u^{(n+1)}(t)\|_{B^{s}_{p, r}} \leq
&e^{C\int_0^t\|u^{(n)}(\tau)\|^2_{B^{s}_{p,
r}}d\tau}\|u_0\|_{B^{s}_{p,
r}}\\&+\frac{C}{\sqrt{2}}\int_0^t
e^{C\int^{t}_{\tau}\|u^{(n)}(\tau')\|^2_{B^{s}_{p,
r}}d\tau^{\prime}} \left(\|u^{(n)}(\tau)\|^3_{B^{s}_{p,
r}}+\|u^{(n)}(\tau)\|_{B^{s}_{p, r}}\right)\;d\tau.
\end{split}
\end{equation}

Let us choose a $T>0$ such that
$$ T\le\min\left\{\frac1{8C\|u_0\|^2_{B^s_{p, r}}},\frac{3(\sqrt{2}-1)}{4C}\right\}, $$
and suppose by induction that for all $ t \in [0, T]$
\begin{equation}\label{3.5}
\|u^{(n)}(t)\|_{B^s_{p, r}} \leq \frac{\sqrt{2}\|u_0\|_{B^s_{p,
r}}}{\left(1-8C\|u_0\|^2_{B^s_{p, r}}t\right)^{1/2}}.
\end{equation}
Indeed, one obtains from \eqref{3.5} that for any $0\leq\tau\leq t$
\begin{align*}
C\int^t_{\tau}&\|u^{(n)}(\tau')\|^2_{B^{s}_{p, r}}\;d\tau' \leq
C\int^t_{\tau}\frac{2\|u_0\|^2_{B^s_{p, r}}}{1-8C\|u_0\|^2_{B^s_{p,
r}}\tau'}\;d\tau'\\&=\frac1{4}\ln(1-8C\|u_0\|^2_{B^s_{p,
r}}\tau)-\frac1{4}\ln(1-8C\|u_0\|^2_{B^s_{p, r}}t).
\end{align*}
And then inserting the above inequality and \eqref{3.5} into
\eqref{3.4} leads to
\begin{align*}
 \|u^{(n+1)}(t)\|_{B^{s}_{p, r}}
\leq&\frac{\|u_0\|_{B^s_{p, r}}}{\sqrt[4]{1-8C\|u_0\|^2_{B^s_{p,
r}}t}}+
 \frac{C}{\sqrt{2}\sqrt[4]{1-8C\|u_0\|^2_{B^s_{p, r}}t}}\\&\times\int_0^t
\left(\frac{2\sqrt{2}\|u_0\|^3_{B^s_{p,
r}}}{(1-8C\|u_0\|^2_{B^s_{p,r}}\tau)^{\frac{5}{4}}}+\frac{\sqrt{2}\|u_0\|_{B^s_{p,
r}}}{(1-8C\|u_0\|^2_{B^s_{p, r}}\tau)^{\frac{1}{4}}}\right)\;d\tau\\
\leq&\frac{\|u_0\|_{B^s_{p, r}}}{\sqrt{1-8C\|u_0\|^2_{B^s_{p,r}}t}}
+\frac{1 -\sqrt[4]{(1-8C\|u_0\|^2_{B^s_{p,
r}}t)^3}}{6\|u_0\|_{B^s_{p,r}}\sqrt[4]{1-8C\|u_0\|^2_{B^s_{p,r}}t}},
\end{align*}
which implies
\begin{align*}
\|u^{(n+1)}(t)\|_{B^s_{p, r}} \leq \frac{\sqrt{2}\|u_0\|_{B^s_{p,
r}}}{\left(1-8C\|u_0\|^2_{B^s_{p, r}}t\right)^{1/2}}.
\end{align*}
Therefore, $(u^{(n)})_{n \in \mathbb{N}}$ is uniformly bounded in
$C([0, T]; B^{s}_{p, r}).$

On the other hand, using the Moser-type estimates (see Lemma
\ref{l2.4} (ii)), one finds that
\begin{align*}
\|[(u^{(n)})^2-(u^{(n)}_x)^2] \pa_x m^{(n+1)}\|_{B^{s-3}_{p,
r}}&\leq C\|m^{(n+1)}\|_{B^{s-2}_{p,
r}}\left(\|u^{(n)}\|^2_{B^{s}_{p, r}}+\|u^{(n)}_x\|^2_{B^{s-1}_{p,
r}}\right)\\&\leq C\|u^{(n+1)}\|_{B^{s}_{p,
r}}\|u^{(n)}\|^2_{B^{s}_{p, r}},
\end{align*}
and
\begin{equation*}
\|u^{(n)}_x (m^{(n)})^2\|_{B^{s-3}_{p, r}}\leq C\|m^{(n)}\|^2_{B^{s-2}_{p,
r}}\|u^{(n)}\|_{B^{s}_{p, r}}\leq C\|u^{(n)}\|^3_{B^{s}_{p, r}}.
\end{equation*}
Hence, using the equation $(T_n),$  we have
\begin{equation*}
\pa_t u^{(n+1)} \in C([0, T]; B^{s-1}_{p, r})
\end{equation*}
uniformly bounded, which yields that the sequence $(u^{(n)})_{n \in
\mathbb{N}}$ is uniformly bounded in $E^{s}_{p, r}(T).$

Next we are going to show that
\begin{equation*}
(u^{(n)})_{n\in \mathbb{N}} \; \mbox{ is a Cauchy sequence in } \;
C([0, T]; B^{s-1}_{p, r}).
\end{equation*}
In fact, according to (\ref{3.3}), we obtain that, for all $n, \, \ell
\in \mathbb{N},$
\begin{equation}\label{e-d-n-1}
\begin{split}
\left\{\pa_t +\left[(u^{(n+\ell)})^2-(u^{(n+\ell)}_x)^2\right]
\pa_x\right\}&(m^{(n+\ell+1)}-m^{(n+1)})\\
&=g(u^{(n+\ell)},u^{(n)},m^{(n+\ell)},m^{(n)},m^{(n+1)}),
\end{split}
\end{equation}
where
\begin{align*}
&g(u^{(n+\ell)},u^{(n)},m^{(n+\ell)},m^{(n)},m^{(n+1)})\\
&\,=[(u^{(n)}-u^{(n+\ell)})(u^{(n)}+u^{(n+\ell)})-(u^{(n)}_x-u^{(n+\ell)}_x)(u^{(n)}_x+u^{(n+\ell)}_x)]\pa_{x}m^{(n+1)}\\
&\,-2u^{(n+\ell)}_x(m^{(n+\ell)}-m^{(n)})(m^{(n)}+m^{(n+\ell)})+2(u^{(n)}_x-u^{(n+\ell)}_x)(m^{(n)})^2+\gamma(u^{(n)}_x-u^{(n+\ell)}_x).
\end{align*}
We rewrite \eqref{e-d-n-1} as the equation in terms of $(u^{(n+\ell+1)}-u^{(n+1)})$
\begin{equation*}
\begin{split}
\left\{\pa_t +\left[(u^{(n+\ell)})^2-(u^{(n+\ell)}_x)^2\right]
\pa_x\right\}&(1-\partial_x^2)(u^{(n+\ell+1)}-u^{(n+1)})\\
&=g(u^{(n+\ell)},u^{(n)},m^{(n+\ell)},m^{(n)},m^{(n+1)}),
\end{split}
\end{equation*}
which is equivalent to
\begin{equation}\label{e-d-n-2}
\begin{split}
(1-\partial_x^2)\bigg\{\left(\pa_t +\left[(u^{(n+\ell)})^2-(u^{(n+\ell)}_x)^2\right]
\pa_x\right)(u^{(n+\ell+1)}-u^{(n+1)})\bigg\}=h^{(n, \ell)}
\end{split}
\end{equation}
with
\begin{equation*}
\begin{split}
h^{(n, \ell)}&=2 \partial_x \left[(u^{(n+\ell)})^2-(u^{(n+\ell)}_x)^2\right] \partial_x^2(u^{(n+\ell+1)}-u^{(n+1)})\\
&\qquad +\partial_x^2 \left[(u^{(n+\ell)})^2-(u^{(n+\ell)}_x)^2\right] \partial_x(u^{(n+\ell+1)}-u^{(n+1)})\\
&\qquad + g(u^{(n+\ell)},u^{(n)},m^{(n+\ell)},m^{(n)},m^{(n+1)}).
\end{split}
\end{equation*}
Applying the operator $(1-\partial_x^2)^{-1}$ to \eqref{e-d-n-2} gives rise to
\begin{equation}\label{e-d-n-3}
\begin{split}
\left\{\pa_t +\left[(u^{(n+\ell)})^2-(u^{(n+\ell)}_x)^2\right]
\pa_x\right\}(u^{(n+\ell+1)}-u^{(n+1)})=(1-\partial_x^2)^{-1} h^{(n, \ell)}
\end{split}
\end{equation}
Thanks to Lemma \ref{l2.2} again, then for every $t \in [0, T],$ we
obtain
 \begin{equation}\label{e-d-n-4}
 \begin{split}
 e&^{-C\int_0^t\| [(u^{(n+\ell)})^2-(u^{(n+\ell)}_x)^2](\tau)\|_{B^{s-1}_{p,
r}}d\tau}\|(u^{(n+\ell+1)}-u^{(n+1)})(t)\|_{B^{s-1}_{p, r}}\\
  &\, \leq
\|u_0^{(n+\ell+1)}-u_0^{(n+1)}\|_{B^{s-1}_{p, r}} +C\int_0^t
e^{-C\int_{0}^{\tau}\|[(u^{(n+\ell)})^2-(u^{(n+\ell)}_x)^2](\tau^{\prime})\|_{B^{s-1}_{p,
r}}d\tau^{\prime}}\|h^{(n, \ell)}\|_{B^{s-3}_{p,
r}}\;d\tau.
\end{split}
\end{equation}
In the case of $s > \max\{2+\frac{1}{p}, \, \frac{5}{2}\}$, one can
deduce from the product law in Besov spaces that
\begin{align*}
\|[(u^{(n)}&-u^{(n+\ell)})(u^{(n)}+u^{(n+\ell)})-(u^{(n)}_x-u^{(n+\ell)}_x)(u^{(n)}_x+u^{(n+\ell)}_x)]\pa_{x}m^{(n+1)}\|_{B^{s-3}_{p, r}}\\
\leq & C\|m^{(n+1)}\|_{B^{s-2}_{p,
r}}(\|u^{(n+\ell)}-u^{(n)}\|_{B^{s-1}_{p,
r}}\|u^{(n+\ell)}+u^{(n)}\|_{B^{s-1}_{p,
r}}\\&\quad+\|u^{(n+\ell)}_x-u^{(n)}_x\|_{B^{s-2}_{p, r}}
\|u^{(n+\ell)}_x+u^{(n)}_x\|_{B^{s-2}_{p, r}})\\
\leq & C\|u^{(n+\ell)}-u^{(n)}\|_{B^{s-1}_{p, r}}
(\|u^{(n)}\|^2_{B^{s}_{p, r}}+\|u^{(n+1)}\|^2_{B^{s}_{p,
r}}+\|u^{(n+\ell)}\|^2_{B^{s}_{p, r}}),
\end{align*}
\begin{align*}
\|u^{(n+l)}_x&(m^{(n+\ell)}-m^{(n)})(m^{(n)}+m^{(n+\ell)})\|_{B^{s-3}_{p,
r}}\\ \leq &C \|u^{(n+\ell)}\|_{B^{s}_{p,
r}}\|m^{(n+\ell)}-m^{(n)}\|_{B^{s-3}_{p,
r}}\|m^{(n+\ell)}+m^{(n)}\|_{B^{s-2}_{p, r}}\\
\leq &C\|u^{(n+\ell)}-u^{(n)}\|_{B^{s-1}_{p, r}}
(\|u^{(n)}\|^2_{B^{s}_{p, r}}+\|u^{(n+\ell)}\|^2_{B^{s}_{p, r}}),
\end{align*}
 and
\begin{align*}
\|(u^{(n)}_x&-u^{(n+\ell)}_x)(m^{(n)})^2\|_{B^{s-3}_{p, r}}\leq
C\|u^{(n+\ell)}-u^{(n)}\|_{B^{s-2}_{p, r}}\|m^{(n)}\|^2_{B^{s-2}_{p,
r}}\\ \leq &C\|u^{(n+\ell)}-u^{(n)}\|_{B^{s-1}_{p, r}}
\|u^{(n)}\|^2_{B^{s}_{p, r}}.
\end{align*}
From this, one finds that
\begin{align*}
\|g&(u^{(n+\ell)},u^{(n)},m^{(n+\ell)},m^{(n)},m^{(n+1)})\|_{B^{s-3}_{p,
r}}\\
& \leq C\|u^{(n+\ell)}-u^{(n)}\|_{B^{s-1}_{p, r}}
\left(\|u^{(n)}\|^2_{B^{s}_{p, r}}+\|u^{(n+1)}\|^2_{B^{s}_{p,
r}}+\|u^{(n+\ell)}\|^2_{B^{s}_{p, r}}+|\gamma|\right).
\end{align*}
Similarly, we may check that
\begin{equation*}
\begin{split}
&\|2 \partial_x \left[(u^{(n+\ell)})^2-(u^{(n+\ell)}_x)^2\right] \partial_x^2(u^{(n+\ell+1)}-u^{(n+1)})\|_{B^{s-3}_{p,
r}} \\
& \leq C\|u^{(n+\ell+1)}-u^{(n+1)}\|_{B^{s-1}_{p, r}}
\|u^{(n+\ell)}\|^2_{B^{s}_{p, r}}
\end{split}
\end{equation*}
and
\begin{equation*}
\begin{split}
&\|\partial_x^2 \left[(u^{(n+\ell)})^2-(u^{(n+\ell)}_x)^2\right] \partial_x(u^{(n+\ell+1)}-u^{(n+1)})\|_{B^{s-3}_{p,
r}}  \\
& \leq C\|u^{(n+\ell+1)}-u^{(n+1)}\|_{B^{s-1}_{p, r}}
\|u^{(n+\ell)}\|^2_{B^{s}_{p, r}}.
\end{split}
\end{equation*}
Hence, we obtain
\begin{equation*}
\begin{split}
&\|h^{(n, \ell)}\|_{B^{s-3}_{p, r}} \leq  C\|u^{(n+\ell+1)}-u^{(n+1)}\|_{B^{s-1}_{p, r}}
\|u^{(n+\ell)}\|^2_{B^{s}_{p, r}}\\
&+C\|u^{(n+\ell)}-u^{(n)}\|_{B^{s-1}_{p, r}}
(\|u^{(n+\ell)}\|^2_{B^{s}_{p, r}}+\|u^{(n+1)}\|^2_{B^{s}_{p, r}}+\|u^{(n)}\|^2_{B^{s}_{p, r}}+|\gamma|).
\end{split}\end{equation*}
Therefore, we obtain
\begin{equation}\label{e-d-n-5}
\begin{split}
e&^{-C\int_0^t\|
[(u^{(n+\ell)})^2-(u^{(n+\ell)}_x)^2](\tau)\|_{B^{s-1}_{p,
r}}d\tau}\|(u^{(n+\ell+1)}-u^{(n+1)})(t)\|_{B^{s-1}_{p, r}}\\
&\leq
\|u_0^{(n+\ell+1)}-u_0^{(n+1)}\|_{B^{s-1}_{p, r}} \\
  & +C\int_0^t
e^{-C\int_{0}^{\tau}\|[(u^{(n+\ell)})^2-(u^{(n+\ell)}_x)^2](\tau^{\prime})\|_{B^{s-1}_{p,
r}}d\tau^{\prime}}\|(u^{(n+\ell)}-u^{(n)})(\tau)\|_{B^{s-1}_{p,
r}} \\
&\times\left(\|u^{(n)}(\tau)\|^2_{B^{s}_{p,
r}}+\|u^{(n+\ell)}(\tau)\|^2_{B^{s}_{p,
r}}+\|u^{(n+1)}(\tau)\|^2_{B^{s}_{p, r}}+|\gamma|\right)\;d\tau\\
  & +C\int_0^t
e^{-C\int_{0}^{\tau}\|[(u^{(n+\ell)})^2-(u^{(n+\ell)}_x)^2](\tau^{\prime})\|_{B^{s-1}_{p,
r}}d\tau^{\prime}}\|(u^{(n+\ell+1)}-u^{(n+1)})(\tau)\|_{B^{s-1}_{p,
r}}\|u^{(n+\ell)}(\tau)\|^2_{B^{s}_{p,
r}}\;d\tau.
\end{split}\end{equation}
Since $(u^{(n)})_{n \in \mathbb{N}}$ is uniformly bounded in
$E^{s}_{p,r}(T)$ and
\begin{equation*}
u_0^{(n+\ell+1)}-u_0^{(n+1)}=S_{n+\ell+1}u_0-S_{n+1}u_0=\sum_{q=n+1}^{n+\ell}\Delta_q
u_0,
 \end{equation*}
then there exists a constant $C_T$ independent of $n$ and $\ell$ such
that for all $t \in [0, T]$
 \begin{equation*}
 \begin{split}
&\|(u^{(n+\ell+1)}-u^{(n+1)})(t)\|_{B^{s-1}_{p, r}}
 \leq  C_T \left (2^{-n}+\int_0^{t}\|(u^{(n+\ell)}-u^{(n)})(\tau)\|_{B^{s-1}_{p, r}}d\tau \right).
 \end{split}
\end{equation*}
Arguing by induction with respect to the index $n$, one can easily prove that
\begin{equation*}
 \begin{split}
\|u^{(n+\ell+1)}-u^{(n+1)}\|_{L^{\infty}_T(B^{s-1}_{p, r})}
 \leq \frac{(T C_T)^{n+1}}{(n+1)!}\|u^{(\ell)}\|_{L^{\infty}_T(B^{s}_{p, r})}+
 C_T \sum_{k=0}^{n}2^{-(n-k)}\frac{(T C_T)^{k}}{k!}.
 \end{split}
 \end{equation*}
Similarly $\|u^{(\ell)}\|_{L^{\infty}_T(B^{s}_{p, r})}$ can be bounded
independently of $\ell$, we conclude that there exist some new constant
$C_T^{\prime}$ independent of $n$ and $\ell$ such that
\begin{equation*}
\|u^{(n+\ell+1)}-u^{(n+1)}\|_{L^{\infty}_T(B^{s-1}_{p, r})} \leq 2^{-n}
C_T^{\prime}.
\end{equation*}
 Hence $(u^{(n)})_{n \in \mathbb{N}}$ is a Cauchy sequence in
$C([0, T]; B^{s-1}_{p, r}).$
\end{proof}

\begin{proof} [{\it Proof of Theorem \ref{t3.1}}.] Thanks to Lemma
\ref{l3.1}, we obtain that $(u^{(n)})_{n \in \mathbb{N}}$ is a
Cauchy sequence in $C([0, T]; B^{s-1}_{p, r}),$ so it converges to
some function $u \in C([0, T]; B^{s-1}_{p, r}).$ We now have to
check that $u$ belongs to $E^s_{ p,r}(T)$  and solves the Cauchy
problem (\ref{e1.1}). Since $(u^{(n)})_{n \in \mathbb{N}}$ is
uniformly bounded in $L^{\infty}([0, T];B^s_{ p,r})$ according to
Lemma \ref{l3.1}, the Fatou property for the Besov spaces
(Proposition \ref{p2.2} iv)) guarantees that $u $ also belongs to
$L^{\infty}([0, T];B^s_{ p,r}).$

On the other hand, as $(u^{(n)})_{n \in \mathbb{N}}$ converges to
$u$ in $C([0, T]; B^{s-1}_{p, r}) ,$ an interpolation argument
ensures that the convergence holds in $C([0, T]; B^{s^{\prime}}_{p,
r}) ,$ for any $s^{\prime} <s$. It is then easy to pass to the limit
in the equation $(T_n)$ and to conclude that $u$ is indeed a
solution to the Cauchy  problem (\ref{e1.1}). Thanks to the fact
that $u$ belongs to $L^{\infty}([0, T];B^s_{ p,r}),$ the right-hand
side of the equation
\begin{equation*}
 \pa_t m +(u^2-u^2_x) \pa_xm=-2u_xm^2-\gamma u_x
\end{equation*}
belongs to $L^{\infty}([0, T];B^{s-2}_{ p,r}).$  In particular, for
the case $r < \infty,$ Lemma \ref{l2.3} implies that $u \in C([0,
T]; B^{s^{\prime}}_{p, r}) $ for any $s^{\prime} < s.$ Finally,
using the equation again, we see that $\pa_t u \,
 \in C([0, T]; B^{s-1}_{p, r})$ if $r < \infty,$ and in
$L^{\infty}([0, T];B^{s-1}_{p,r})$ otherwise.
 Moreover, a standard use
of a sequence of viscosity approximate solutions
$(u_{\epsilon})_{\epsilon
>0}$ for the Cauchy  problem (\ref{e1.1}) which converges uniformly in
$$C([0, T]; B^{s}_{p, r}) \cap C^{1}([0, T];
B^{s-1}_{p, r}) $$ leads to the continuity of the solution $u$ in
$E^s_{ p,r}(T) .$
\end{proof}

\subsection{Critical case}
Attention is now restricted the critical case in the local well-posedness.
\begin{thm}\label{t3.2}
Suppose that the initial data $u_0(x)\in B^{\frac{5}{2}}_{2,1}$.
Then there exists a maximal $T=T(u_0)>0$ and a unique solution $u(t,
x)$ to the Cauchy problem (\ref{e1.1}) such that
$$u=u(\cdot,u_0)\in C([0,T];B^{\frac{5}{2}}_{2,1})\cap C^1([0,T];B^{\frac{3}{2}}_{2,1}).$$
Moreover, the solution depends continuously on the initial data,
i.e. the mapping
$$u_0\mapsto u(\cdot,u_0):B^{\frac{5}{2}}_{2,1}\mapsto C([0,T];B^{\frac{5}{2}}_{2,1})
\cap C^1([0,T];B^{\frac{3}{2}}_{2,1})$$ is continuous.
\end{thm}

\begin{rmk}\label{rmk3.5}
Note that the equation in \eqref{e1.1} with regard to $m$ is a
transport form, that is,
$$m_t+(u^2-u^2_x)m_x=-2u_xm^2-\gamma u_x.$$
Roughly speaking, in order to propagate the regularity of the
solution $m$ to the Cauchy problem (\ref{e1.1}) in terms of its
initial data $m_0$, the ``coefficient'' $u^2-u^2_x$ of $m_x$ needs
to satisfy the Lipschitz condition. Toward this purpose, it suffices
to guarantee $u$ belonging to $W^{2, \infty}$, the space of bounded
functions with bounded first and second derivatives, which
satisfies the embedding properties $
B^{\frac{5}{2}}_{2,1}\hookrightarrow B^2_{\infty,1}\hookrightarrow
W^{2, \infty} \hookrightarrow B^2_{\infty,\infty}$. From this, we
call $s=\frac{5}{2}$ the critical regularity index in terms of $u$
for the well-posedness of the initial value problem \eqref{e1.1} in
the following sense:
\begin{equation*}H^s\hookrightarrow B^{\frac{5}{2}}_{2,1}\hookrightarrow H^{\frac{5}{2}}\hookrightarrow
 B^{\frac{5}{2}}_{2,\infty}\hookrightarrow H^{s'} \quad \mbox{for all} \quad s'<\frac{5}{2}<s.
 \end{equation*}
\end{rmk}

\begin{rmk}\label{rmk3.6} Similar to the result of the Camassa-Holm equation presented by Danchin in \cite{Dan1},
 using the estimates in the proofs of Theorem \ref{t3.1}-\ref{t3.2}, we may demonstrate the well-posedness of
the equation \eqref{e1.1} with the initial data $u_0$ belonging to
the critical space $B^{\frac{5}{2}}_{2, \infty} \cap W^{2, \infty}$.
We leave the details to the readers.
\end{rmk}

\begin{proof} [Proof of Theorem \ref{t3.2}]  Theorem \ref{t3.2} will be divided into the following three lemmas.
\end{proof}

We first present the existence of the solution.
\begin{lem}\label{l4.1}
Assume that $u_0\in B^{\frac{5}{2}}_{2,1}$. Then there
exists a time $T>0$ such that the Cauchy problem (\ref{e1.1}) has a
solution $u\in C([0,T];B^{\frac{5}{2}}_{2,1}) \cap
C^1([0,T];B^{\frac{3}{2}}_{2,1})$.
\end{lem}
\begin{proof}
On account of $u_0\in B^{\frac{5}{2}}_{2,1}$, the transport
theory (see Lemma \ref{l2.2}) can be applied. Similar to the case
$u_0(x)\in B^s_{p,r},\, s> \max\{2+\frac{1}{p}, \, \frac{5}{2}\}$ ,
we can establish this lemma. The proof of the lemma is therefore omitted without details.
\end{proof}
We are now in a position to establish estimates in
$L^{\infty}(0,T;B^{\frac{3}{2}}_{2,1})$ for the difference of two
solutions of the Cauchy problem \eqref{e1.1} belonging to
$L^{\infty}([0,T];B^{\frac{5}{2}}_{2,1})\cap
C([0,T];B^{\frac{3}{2}}_{2,1})$. Uniqueness is a corollary of the
following result.
\begin{lem}\label{l4.2}
Suppose that $u_0$ (resp. $v_0$) $\in B^{\frac{5}{2}}_{2,1}$ such
that $u$ (resp. $v$) $\in
L^{\infty}([0,T];B^{\frac{5}{2}}_{2,1})\cap
C([0,T];B^{\frac{3}{2}}_{2,1})$ is a solution to the Cauchy problem
(\ref{e1.1}) with initial data $u_0$ (resp. $v_0$). Let $w= u-v$ and
$w_0=u_0-v_0$.  Then for every $t \in [0, T]:$
\begin{align}\label{4.7}
\|w(t)\|_{B^{\frac{3}{2}}_{2,1}} \leq
\|w(0)\|_{B^{\frac{3}{2}}_{2,1}} \exp \left
\{C\int_{0}^{t}\left(\|u(\tau)\|^2_{B^{\frac{5}{2}}_{2,1}}+\|v(\tau)\|^2_{B^{\frac{5}{2}}_{2,1}}+|\gamma|\right)\;d\tau \right\} .
\end{align}
\end{lem}
\begin{proof}
 Thanks to the formulation \eqref{eqns-2}, we see that $w$ solves the linear
equation
\begin{equation}\label{eqns-d-1}
\begin{split}
&\partial_t w+\left[u^2-\frac{1}{3}(u_x^2+u_xv_x+v^2_x)\right]\pa_x
w+(1-\partial_x^2)^{-1}\bigg(\frac{1}{3}(u_x^2+u_xv_x+v^2_x)\, w_x+ \gamma \, w_x\bigg)\\
&+(u+v)v_x \,
w+\partial_x(1-\partial_x^2)^{-1}\bigg(\frac{2}{3}(u^2+uv+v^2)\, w+
u(u_x+v_x) \,w_x +v_x^2\, w\bigg)=0.
\end{split}
\end{equation}
Consider that $u_0$ (resp. $v_0$) $\in B^{\frac{5}{2}}_{2,1}$ such
that $u$ (resp. $v$) $\in
L^{\infty}([0,T];B^{\frac{5}{2}}_{2,1})\cap
C([0,T];B^{\frac{3}{2}}_{2,1})$, by virtue of the transport theory
in Lemma \ref{l2.2}, the following inequality holds true:
\begin{equation}\label{eqns-d-2}
\begin{split}
\|w(t)\|_{B^{\frac3{2}}_{2,1}}\leq \|w(0)\|_{B^{\frac3{2}}_{2,1}}&+
C\int_0 ^t\left \|u^2-\frac{1}{3}(u_x^2+u_xv_x+v^2_x)\right
\|_{B^{\frac3{2}}_{2, 1}}(\tau)\,\|w(\tau)\|_{B^{\frac3{2}}_{2,1}}
d\tau\\
&+ \int_0^t\|f(u,v, u_x, v_x, w, w_x)(\tau)\|_{B^{\frac3{2}}_{2,1}}\,d\tau,
\end{split}
\end{equation}
where
\begin{equation*}
\begin{split}
&f(u,v, u_x, v_x, w, w_x)=(1-\partial_x^2)^{-1}\bigg(\frac{1}{3}(u_x^2+u_xv_x+v_x^2)\, w_x+ \gamma \, w_x\bigg)\\
&+(u+v)v_x \,
w+\partial_x(1-\partial_x^2)^{-1}\bigg(\frac{2}{3}(u^2+uv+v^2)\, w+
u(u_x+v_x) \,w_x +v^2_x\, w\bigg).
\end{split}
\end{equation*}
Applying the product law in the Besov spaces, we have
\begin{equation*}
\begin{split}
\left\|u^2-\frac{1}{3}(u_x^2+u_xv_x+v^2_x)\right\|_{B^{\frac3{2}}_{2,
1}} \leq C (\|u\|_{B^{\frac5{2}}_{2, 1}}^2+\|v\|_{B^{\frac5{2}}_{2,
1}}^2).
\end{split}
\end{equation*}
Similarly, one gets
\begin{equation*}
\begin{split}
\left\|(1-\partial_x^2)^{-1} \right . & \left .  \bigg(\frac{1}{3}(u_x^2+u_xv_x+v_x^2)\, w_x+ \gamma \, w_x\bigg)\right\|_{B^{\frac{3}{2}}_{2,1}}\\
&\leq C\left\|\frac{1}{3}(u_x^2+u_xv_x+v_x^2)\, w_x+ \gamma \, w_x\right\|_{B^{\frac{1}{2}}_{2,1}}\\
&\leq
C\bigg(\|u\|_{B^{\frac{3}{2}}_{2,1}}^2+\|u\|_{B^{\frac{3}{2}}_{2,1}}\|v\|_{B^{\frac{3}{2}}_{2,1}}
+\|v\|_{B^{\frac{3}{2}}_{2,1}}^2+|\gamma|\bigg)\,
\|w\|_{B^{\frac{3}{2}}_{2,1}},
\end{split}
\end{equation*}
$$
\|(u+v)v_x \, w\|_{B^{\frac{3}{2}}_{2,1}} \leq C(\|u\|_{B^{\frac{3}{2}}_{2,1}}+\|v\|_{B^{\frac{3}{2}}_{2,1}})
\|v \|_{B^{\frac{5}{2}}_{2,1}}\| w\|_{B^{\frac{3}{2}}_{2,1}},
$$
and
\begin{equation*}
\begin{split}
\left\|\partial_x \right .  & \left . (1-\partial_x^2)^{-1}  \bigg(\frac{2}{3}(u^2+uv+v^2)\, w+ u(u_x+v_x) \,w_x
 +v^2_x\, w\bigg)\right\|_{B^{\frac{3}{2}}_{2,1}}\\
&\leq C\left\|\frac{2}{3}(u^2+uv+v^2)\, w+ u(u_x+v_x) \,w_x +v^2_x\, w\right\|_{B^{\frac{1}{2}}_{2,1}}\\
&\leq
C\bigg(\|u\|_{B^{\frac{1}{2}}_{2,1}}^2+\|v\|_{B^{\frac{1}{2}}_{2,1}}^2
+\|u\|_{B^{\frac{1}{2}}_{2,1}}(\|u\|_{B^{\frac{3}{2}}_{2,1}}+\|v\|_{B^{\frac{3}{2}}_{2,1}})+\|v\|^2_{B^{\frac{3}{2}}_{2,1}}\bigg)\,
 \|w\|_{B^{\frac{3}{2}}_{2,1}},
\end{split}
\end{equation*}
which leads to
\begin{equation*}
\begin{split}
&\|f(u,v, u_x, v_x, w, w_x)\|_{B^{\frac{3}{2}}_{2,1}}\leq
C\bigg(\|u\|_{B^{\frac{5}{2}}_{2,1}}^2+\|v\|_{B^{\frac{5}{2}}_{2,1}}^2+|\gamma|\bigg)\,
\|w\|_{B^{\frac{3}{2}}_{2,1}}.
\end{split}
\end{equation*}
Hence, we obtain from \eqref{eqns-d-2} that
\begin{align*}
\|w(t)\|_{B^{\frac{3}{2}}_{2,1}} \leq
\|w(0)\|_{B^{\frac{3}{2}}_{2,1}}
+C\int_{0}^{t}\bigg(\|u(\tau)\|_{B^{\frac{5}{2}}_{2,1}}^2+\|v(\tau)\|_{B^{\frac{5}{2}}_{2,1}}^2
+|\gamma|\bigg) \|w(\tau)\|_{B^{\frac{3}{2}}_{2,1}}\, d\tau.
\end{align*}
Therefore, due to the Gronwall inequality, we deduce that
\begin{align*}
\|w(t)\|_{B^{\frac{3}{2}}_{2,1}} \leq
\|w(0)\|_{B^{\frac{3}{2}}_{2,1}} \exp \left
\{C\int_{0}^{t}\bigg(\|u(\tau)\|_{B^{\frac{5}{2}}_{2,1}}^2+\|v(\tau)\|_{B^{\frac{5}{2}}_{2,1}}^2
+|\gamma|\bigg)\;d\tau \right\} .
\end{align*}
which is the desired result.
\end{proof}
At last, we are going to verify the continuity of the solution with
regard to initial data in $B^{\frac{5}{2}}_{2,1}$.
\begin{lem}\label{l4.3}
For any $u_0\in B^{\frac{5}{2}}_{2,1}$, there exists a time $T>0$
and a neighborhood $V$ of $u_0$ in $B^{\frac{5}{2}}_{2,1}$ such that
for any $v\in V$, which is the solution of the Cauchy problem
(\ref{e1.1}) with the initial data $v_0$, the map
\begin{align*}
\Phi: v_0\rightarrow v(\cdot,v_0): V\subset
B^{\frac{5}{2}}_{2,1}\rightarrow C([0,T];B^{\frac{5}{2}}_{2,1}) \cap
C^1([0,T];B^{\frac{3}{2}}_{2,1})
\end{align*}
is continuous.
\end{lem}
Motivated by \cite{Dan2}, Lemma \ref{l4.3} can be established by applying Lemma
\ref{l4.2} and a continuity result Lemma \ref{l2.5} for the linear
transport equations.
\begin{proof}[Proof of Lemma \ref{l4.3}]
 We first prove the continuity of the map $\Phi$ in $C([0,T]; B^{\frac{3}{2}}_{2,1})$. Let us
fix a $u_0\in B^{\frac{5}{2}}_{2,1}$ and a $\delta>0$. We claim that
there exists a $T>0$ and $M>0$ such that for any $\tilde{u}_0\in
B^{\frac{5}{2}}_{2,1}$ with
$\|\tilde{u}_0-u_0\|_{B^{\frac{5}{2}}_{2,1}}\leq\delta$, the solution
$\tilde{u}=\Phi(\tilde{u}_0)$ of the Cauchy problem (\ref{e1.1}) associated to
$\tilde{u}_0$ belongs to $C([0,T];B^{\frac{5}{2}}_{2,1})$ and satisfies
 $ \displaystyle \|\tilde{u}\|_{L^{\infty}(0,T;B^{\frac{5}{2}}_{2,1})}\leq M.$
Indeed, from the proof of the local well-posedness we know that when
we fix a $T>0$ such that
$$T\le\min\left\{\frac1{8C\|\tilde{u}_0\|^2_{B^{\frac{5}{2}}_{2,1}}},\frac{3(\sqrt{2}-1)}{4C}\right\},$$
then from \eqref{3.5} similarly we deduce that
\begin{align}\label{e4.5}
\|\tilde{u}(t)\|_{B^{\frac{5}{2}}_{2,1}} \leq
\frac{\sqrt{2}\|\tilde{u}_0\|_{B^{\frac{5}{2}}_{2,1}}}{\left(1-8C\|\tilde{u}_0\|^2_{B^{\frac{5}{2}}_{2,1}}t\right)^{1/2}}\qquad\text{for
all}\qquad t\in[0,T].
\end{align}

Since $\|\tilde{u}_0-u_0\|_{B^{\frac{5}{2}}_{2,1}}\leq\delta$, we get
$\|\tilde{u}_0\|_{B^{\frac{5}{2}}_{2,1}}\leq\|u_0\|_{B^{\frac{5}{2}}_{2,1}}+\delta$.
One can choose some suitable constant $C$, such that
$$T=\frac{3}{32C\left(\|u_0\|_{B^{\frac{5}{2}}_{2,1}}+\delta+1\right)^2}\qquad\text{and}\qquad
 M=2\sqrt{2}(\|u_0\|_{B^{\frac{5}{2}}_{2,1}}+\delta).$$
 Now combining the above uniform bounds with Lemma
\ref{l4.2}, we infer that
$$\|\Phi(\tilde{u}_0)-\Phi(u_0)\|_{L^{\infty}(0,T;B^{\frac{3}{2}}_{2,1})}
\leq \delta e^{C(2M^2+|\gamma|)T}.$$ In view of this inequality, we
know that $\Phi$ is the H\"{o}lder continuous from
$B^{\frac{5}{2}}_{2,1}$ into $C([0,T];B^{\frac{3}{2}}_{2,1})$.

Next, we present the continuity of the map $\Phi$ in $C([0,T];B^{\frac{5}{2}}_{2,1})$. Let
$u^{(\infty)}_0\in B^{\frac{5}{2}}_{2,1}$ and
$(u_0^{(n)})_{n\in\mathbb{N}}$ tend to $u^{(\infty)}_0$ in
$B^{\frac{3}{2}}_{2,1}$. Denote by $u^{(n)}$ the solution
corresponding to datum $u_0^{(n)}$. According to the above argument, one can
find $T,\,M>0$ independent of $n$ such that for all $n\in\mathbb{N}$, $u^{(n)}$ is
defined on $[0,T]$ and
\begin{equation}\label{diff-2}
\sup\limits_{n\in\mathbb{N}}\|u^{(n)}\|_{L^{\infty}_T(B^{\frac{5}{2}}_{2,1})}\leq
M.
\end{equation} Thanks to the first step, proving that $u^{(n)}$ tends to
$u^{(\infty)}$ in $C([0,T];B^{\frac{5}{2}}_{2,1})$ amounts to
proving that $m^{(n)}=u^{(n)}-u^{(n)}_{xx}$ tends to $
m^{(\infty)}=u^{(\infty)}-u^{(\infty)}_{xx}$ in
$C([0,T];B^{\frac{1}{2}}_{2,1})$.

Note that $m^{(n)}$ solves the following linear transport equation:
\begin{equation*}
\begin{cases}
\left\{\pa_t +\left[(u^{(n)})^2-(u^{(n)}_x)^2\right]
\pa_x\right\}m^{(n)}=f^{(n)},\\
m^{(n)}|_{t=0} =m_{0}^{(n)}(x)=u_0^{(n)}-u^{(n)}_{0xx},
\end{cases}
\end{equation*}
where $f^{(n)}=-2u_x^{(n)}[m^{(n)}]^2-\gamma u_x^{(n)}.$ Following
the Kato theory \cite{kato}, we decompose $m^{(n)}$ into
$m^{(n)}=z^{(n)}+w^{(n)}$ with
\begin{equation*}
\begin{cases}
\left\{\pa_t +\left[(u^{(n)})^2-(u^{(n)}_x)^2\right]
\pa_x\right\}z^{(n)}=f^{(n)}-f^{(\infty)},\\
z^{(n)}|_{t=0} =m_{0}^{(n)}(x)-m_{0}^{(\infty)}(x).
\end{cases}
\end{equation*}
and
\begin{equation}\label{diff-eqns-1}
\begin{cases}
\left\{\pa_t +\left[(u^{(n)})^2-(u^{(n)}_x)^2\right]
\pa_x\right\}w^{(n)}=f^{(\infty)},\\
w^{(n)}|_{t=0} =m_{0}^{(\infty)}(x).
\end{cases}
\end{equation}
Using Lemma \ref{l2.3} and the product law in the Besov spaces, one may check that
\begin{equation}\label{diff-1}
\begin{split}
&\|z^{(n)}(t)\|_{B^{\frac{1}{2}}_{2,1}}\leq
\exp\{C\int_0^t\|[(u^{(n)})^2-(u_x^{(n)})^2](\tau)\|_{B^{\frac{3}{2}}_{2,
1}} d\tau\}\\
&\qquad\qquad \qquad\qquad\times\bigg(\|m_0^{(n)}-m_0^{(\infty)}\|_{B^{\frac{1}{2}}_{2,1}}
+C\int_0^t\|(f^{(n)}-f^{(\infty)})(\tau)\|_{B^{\frac1{2}}_{2,
1}}\, d\tau\bigg)\\
&\leq \exp\{C\int_0^t\|u^{(n)}(\tau)\|_{B^{\frac{5}{2}}_{2,
1}}^2\, d\tau\}\bigg[\|m_0^{(n)}-m_0^{(\infty)}\|_{B^{\frac{1}{2}}_{2,1}}+C\int_0^t\left(\|u^{(n)}\|^2_{B^{\frac{5}{2}}_{2,
1}}+\|u^{(\infty)}\|^2_{B^{\frac{5}{2}}_{2,
1}}+|\gamma|\right)\\
&\qquad \qquad\qquad\qquad\times\left(\|(m^{(n)}-m^{(\infty)})(\tau)\|_{B^{\frac{1}{2}}_{2,
1}}+\|(u^{(n)}-u^{(\infty)})(\tau)\|_{B^{\frac{3}{2}}_{2, 1}}\right)
d\tau\bigg].
\end{split}
\end{equation}
On the other hand, since the sequence
$(u^{(n)})_{n\in\bar{\mathbb{N}}}$ is uniformly bounded in
$C([0,T];B^{\frac{5}{2}}_{2,1})$ and tends to $u^{(\infty)}$ in
$C([0,T];B^{\frac{3}{2}}_{2,1})$, applying Lemma \ref{l2.5} to \eqref{diff-eqns-1} implies that
$w^{(n)}$ tends to $m^{(\infty)}$ in
$C([0,T];B^{\frac{1}{2}}_{2,1})$.

Let $\varepsilon>0$. Thanks to the above result of convergence with
estimates \eqref{diff-2} and \eqref{diff-1}, we deduce that for large enough
$n\in\mathbb{N}$
\begin{align*}
\|(m^{(n)}&-m^{(\infty)})(\tau)\|_{B^{\frac{1}{2}}_{2, 1}}\\  \leq
 & \, \varepsilon + C(M^2+|\gamma|)
e^{C(M^2+|\gamma|)t}\Big\{\|m_0^{(n)}-m_0^{(\infty)}\|_{B^{\frac{1}{2}}_{2,1}}\\&+\int_0^t
\left(\|(m^{(n)}-m^{(\infty)})(\tau)\|_{B^{\frac{1}{2}}_{2,
1}}+\|(u^{(n)}-u^{(\infty)})(\tau)\|_{B^{\frac{3}{2}}_{2,
1}}\right)d\tau\Big\}.
\end{align*}
As $u^{(n)}$ tends to $u^{(\infty)}$ in
$C([0,T];B^{\frac{3}{2}}_{2,1})$, the last term in the above
integral is less than $\varepsilon$ for large $n$. Hence, thanks to
Gronwall's inequality, we get
\begin{align*}
\|(m^{(n)}-m^{(\infty)})(\tau)\|_{L^{\infty}(0,T;B^{\frac{1}{2}}_{2,
1})}\leq
C_{M,T,|\gamma|}\left(\varepsilon+\|m_0^{(n)}-m_0^{(\infty)}\|_{B^{\frac{1}{2}}_{2,1}}\right)
\end{align*}
for some constant $C_{M,T,|\gamma|}$ depending only on $M$,
$|\gamma|$ and $T$,  which completes the continuity of the map $\Phi$ in $C([0,T];B^{\frac{5}{2}}_{2,1})$.

Finally, applying $\partial_t$ to the equation in (\ref{e1.1}) and using the same argument to the resulting
 equation in terms of $\partial_t u$, we may verify the continuity of the map $\Phi$ in $C^1([0,T];B^{\frac{3}{2}}_{2,1})$.
  This completes the proof of Lemma \ref{l4.3}.
\end{proof}

\renewcommand{\theequation}{\thesection.\arabic{equation}}
\setcounter{equation}{0}

\section{Blow-up scenario and a lower bound of the maximal existence time}

In \cite{gloq}, the authors derived a new wave-breaking mechanism
for solutions to the equation in \eqref{e1.1} with certain initial
profiles. It means that the maximal time of existence of solutions
to the equation in \eqref{e1.1} has some definite upper bounds under
given initial conditions. In this section, we will give a lower
bound of the existence time for this equation.  We  first present
the following Theorem.
\begin{thm}\label{t4.1}(\cite{gloq})
Let $m_0=(1-\partial_x^2)u_0\in H^s(\mathbb{R})$ with $ s>
\frac{1}{2}$.  Let $m$ be the corresponding solution to
\eqref{e1.1}. Assume $ T^{\ast}_{m_0}
> 0 $ is the maximum time of existence. Then
\begin{equation}\label{4.1}
T^{\ast}_{u_0} < \infty \, \, \Rightarrow \, \,
\int_{0}^{T^{\ast}_{m_0}} \|(mu_x)(\tau)\|_{L^{\infty}}\,d\tau =
\infty.
\end{equation}
\end{thm}
\begin{rmk}
There is a little difference between this theorem and the original theorem (Theorem 4.2 in
\cite{gloq}). We recall the result in \cite{gloq} as follows:
\begin{equation}\label{4.1-gloq}
T^{\ast}_{u_0} < \infty \, \, \Rightarrow \, \,
\int_{0}^{T^{\ast}_{m_0}} \|m(\tau)\|_{L^{\infty}}^2\,d\tau =
\infty,
\end{equation}
which, together with the maximum principle to the transport equation \eqref{e1.1} in terms of $m$ applied, implies \eqref{4.1}.
In fact,  applying the maximum principle to the transport equation \eqref{e1.1}, we immediately get
\begin{equation}\label{maximum-1}
\|m(t)\|_{L^{\infty}} \leq \|m_0\|_{L^{\infty}} + C\int_0^{t}(\|(mu_x)(\tau)\|_{L^{\infty}}+|\gamma|)\|m(\tau)\|_{L^{\infty}}\,d\tau,
\end{equation}
where we used the estimate $\|u_x\|_{L^{\infty}} \leq C \|m\|_{L^{\infty}}$ and $\|(u^2-u_x^2)_x\|_{L^{\infty}}=2\|mu_x\|_{L^{\infty}}$. Then,
Gronwall's inequlity applied to \eqref{maximum-1} yields
\begin{equation}\label{maximum-2}
\|m(t)\|_{L^{\infty}} \leq \|m_0\|_{L^{\infty}} \exp\{ C\int_0^{t}(\|(mu_x)(\tau)\|_{L^{\infty}}+|\gamma|)\,d\tau\},
\end{equation}
which along with \eqref{4.1-gloq} gives rise to \eqref{4.1}.
\end{rmk}

In the following, attention is now turned to blow-up issue.  We
first recall a blow-up scenario in \cite{gloq}.
\begin{thm}\label{t4.2}(\cite{gloq})
Let $u_0\in H^s$, $s >5/2$, and $u(t,x)$ be the solution of the
Cauchy  problem (\ref{e1.1}) with life-span $T$. Then $T$ is finite
if and only if
$$\underset{t\uparrow T}{\liminf}\left (\underset{x\in \mathbb{R}}{\inf}
(mu_x(t,x))\right )=-\infty.$$
\end{thm}

We now deduce a lower bound depending only on $\|u_0\|_{W^{2,
\infty}}$ for the maximal time of existence of the solution to
\eqref{e1.1}.

\begin{thm}\label{t4.3}
Assume that $u_0\in H^s$ with $ s>
\frac{5}{2}$. Let  $ T^{\ast}
> 0 $ be the maximum time of existence of the solution $u$ to \eqref{e1.1} with the initial data $u_0$.
If $\gamma\neq0$, then $T^{\ast}$ satisfies
\begin{equation*}
T^{\ast}
\ge\dfrac1{20|\gamma|}\ln\left(1+\dfrac{2|\gamma|}{(2\|u_0\|_{L^{\infty}}+\|\pa_xu_0\|_{L^{\infty}}+2\|\pa^2_xu_0\|_{L^{\infty}})^2}\right).
\end{equation*}
Otherwise, if $\gamma=0$, then
\begin{equation*}
T^{\ast}
\ge\dfrac{1}{2(3\|u_0\|_{L^{\infty}}+3\|\pa_xu_0\|_{L^{\infty}}+\|\pa^2_xu_0\|_{L^{\infty}})^2}.
\end{equation*}
\end{thm}
\begin{proof}
Note that the equation in \eqref{e1.1} is equivalent to the
following equation
\begin{equation}\label{4.2}
u_t+\left(u^2-\frac1{3}u^2_x\right)u_x+\partial_x
G*\left(\frac{2}{3}u^3+uu^2_x\right)+G*\left(\frac1{3}u^3_x+\gamma
u_x\right)=0,
\end{equation}
where $u=G*m=(1-\partial^2_x)^{-1}m$ and $G(x)=\frac1{2}e^{-|x|}$.
Multiplying the above equation by $u^{2n-1}$ and integrating the
results in $x$-variable, in view of H\"{o}lder's inequality, we
obtain
\begin{align*}
&\int_{\mathbb{R}}u^{2n-1}u_t \ dx= \dfrac{1}{2n}\dfrac{\ d}{\ d
t}\|u\|^{2n}_{L^{2n}}
=\|u\|^{2n-1}_{L^{2n}}\frac{\ d}{\ d t}\|u\|_{L^{2n}},\\
&\left |\int_{\mathbb{R}}u^{2n-1}u^3_x\
dx\right|\le\|u_x\|^2_{L^{\infty}}\|u\|^{2n-1}_{L^{2n}}\|u_x\|_{L^{2n}},
\end{align*}
and
\begin{align*}
 & \left|\int_{\mathbb{R}}u^{2n-1}
G\ast\left(\partial_x\left(\frac{2}{3}u^3+uu^2_x\right)+\frac1{3}u^3_x+\gamma
u_x\right)\ d
x\right|\\&\quad\le\|u\|^{2n-1}_{L^{2n}}\left(\frac1{3}\|u_x\|_{L^{\infty}}^2\|u_x\|_{L^{2n}}+
\frac{2}{3}\|u\|_{L^{\infty}}^2\|u\|_{L^{2n}}+\|u_x\|_{L^{\infty}}^2\|u\|_{L^{2n}}+|\gamma|\|u_x\|_{L^{2n}}\right).
\end{align*}
Integrating over $[0,t]$, it follows that
\begin{align*}
\|u(t)\|_{L^{2n}}\le\|u(0)\|_{L^{2n}}&+\int^t_0\|u(\tau)\|_{L^{2n}}\left(\|u_x(\tau)\|_{L^{\infty}}^2+
\frac{2}{3}\|u(\tau)\|_{L^{\infty}}^2\right)\
d\tau\\&\quad+\frac{4}{3}\int^t_0\|u_x(\tau)\|_{L^{2n}}(\|u_x(\tau)\|_{L^{\infty}}^2+|\gamma|)\
d\tau.
\end{align*}
Letting $n$ tend to infinity in the above inequality, we have
\begin{align}\label{4.3}
\|u(t)\|_{L^{\infty}}\le\|u(0)\|_{L^{\infty}}&+\int^t_0\left(\|u(\tau)\|_{L^{\infty}}\|u_x(\tau)\|_{L^{\infty}}^2+
\frac{2}{3}\|u(\tau)\|_{L^{\infty}}^3+|\gamma|\|u_x(\tau)\|_{L^{\infty}}\right)\
d\tau\nonumber\\&\quad+\frac{4}{3}\int^t_0\|u_x(\tau)\|_{L^{\infty}}^3\
d\tau.
\end{align}
Differentiating \eqref{4.2} with respect to $x$, in view of
$(1-\partial^2_x)G*f=f$, we obtain
\begin{equation}\label{4.4}
u_{tx}+u^2u_{xx}+uu_x^2-\frac{2}{3}u^3-u^2_xu_{xx}+
G*\left(\frac{2}{3}u^3+uu^2_x\right)+\partial_xG*\left(\frac1{3}u^3_x+\gamma
u_x\right)=0.
\end{equation}
Multiplying the above equation by $u_x^{2n-1}$ and integrating the
results in $x$ over ${\mathbb R}$, still in view of H\"{o}lder's
inequality, we have
\begin{equation*}
\begin{split}
\int_{\mathbb{R}}u_x^{2n-1}u_{xt} \ dx
&=\|u_x\|^{2n-1}_{L^{2n}}\frac{\ d}{\ d t}\|u_x\|_{L^{2n}},\\
\left|\int_{\mathbb{R}}u_x^{2n-1}u^2u_{xx}\ dx \right|&=
\left|\dfrac{1}{2n}\int_{\mathbb{R}}u^2(u_{x}^{2n})_x\ dx\right|=\left|-\dfrac{1}{2n}\int_{\mathbb{R}}2uu_{x}^{2n+1}\ dx\right|\\
&\le\dfrac{2}{2n}\|u_x\|^2_{L^{\infty}}\|u_x\|^{2n-1}_{L^{2n}}\|u\|_{L^{2n}},\\
\left|\int_{\mathbb{R}}uu_x^{2n+1}\ dx \right|&
\le\|u_x\|^2_{L^{\infty}}\left|\int_{\mathbb{R}}u_x^{2n-1}u\ d
x\right|
\le\|u_x\|^2_{L^{\infty}}\|u_x\|^{2n-1}_{L^{2n}}\|u\|_{L^{2n}},\\
\left|\int_{\mathbb{R}}u^3u_x^{2n-1}\ dx\right|
&\le\|u\|^2_{L^{\infty}}\left|\int_{\mathbb{R}}u_x^{2n-1}u\ d
x\right|
\le\|u\|^2_{L^{\infty}}\|u_x\|^{2n-1}_{L^{2n}}\|u\|_{L^{2n}},
\end{split}
\end{equation*}
and
\begin{align*}
& \left|\int_{\mathbb{R}}u_x^{2n-1}
G\ast\left(\frac{2}{3}u^3+uu^2_x+\pa_x\left(\frac1{3}u^3_x+\gamma
u_x\right)\right)\
dx\right|\\&\quad\le\|u_x\|^{2n-1}_{L^{2n}}\left(\frac1{3}\|u_x\|_{L^{\infty}}^2\|u_x\|_{L^{2n}}+
\frac{2}{3}\|u\|_{L^{\infty}}^2\|u\|_{L^{2n}}+\|u_x\|_{L^{\infty}}^2\|u\|_{L^{2n}}+|\gamma|\|u_x\|_{L^{2n}}\right).
\end{align*}
Integrating over $[0,t]$ and letting $n$ tend to infinity, it
follows that
\begin{align}\label{4.5}
\|u_x(t)\|_{L^{\infty}}\le\|u_x(0)\|_{L^{\infty}}&+\int^t_0\left(3\|u(\tau)\|_{L^{\infty}}\|u_x(\tau)\|_{L^{\infty}}^2+
\frac{4}{3}\|u(\tau)\|_{L^{\infty}}^3+|\gamma|\|u_x\|_{L^{\infty}}\right)\
d\tau\nonumber\\&\quad+\frac1{3}\int^t_0\|u_x(\tau)\|_{L^{\infty}}^3\
d\tau.
\end{align}
Differentiating \eqref{4.4} with respect to $x$, in view of
$(1-\partial^2_x)G*f=f$, we obtain
\begin{align}\label{4.6}
&u_{txx}+u^2u_{xxx}+4uu_xu_{xx}-2u^2u_x+\frac{2}{3}u_x^3-2u_xu_{xx}^2-u^2_xu_{xxx}-\gamma
u_x\nonumber\\&\quad+
\partial_xG*\left(\frac{2}{3}u^3+uu^2_x\right)+G*\left(\frac1{3}u^3_x+\gamma
u_x\right)=0.
\end{align}
Multiplying the above equation by $u_{xx}^{2n-1}$ and integrating
the results in $x$-variable, still in view of H\"{o}lder's
inequality, we have
\begin{equation*}
\begin{split}
&\int_{\mathbb{R}}u_{xx}^{2n-1}u_{xxt} \ dx
=\|u_{xx}\|^{2n-1}_{L^{2n}}\frac{\ d}{\ d t}\|u_{xx}\|_{L^{2n}},
\end{split}
\end{equation*}
\begin{equation*}
\begin{split}
\left|\int_{\mathbb{R}}u_{xx}^{2n-1}u^2u_{xxx}\ dx \right| &=
\left|\dfrac{1}{2n}\int_{\mathbb{R}}u^2(u_{xx}^{2n})_x\ dx\right|=\left|-\dfrac{1}{2n}\int_{\mathbb{R}}2uu_xu_{xx}^{2n}\ dx\right|\\
&\le\dfrac{2}{2n}\|u\|_{L^{\infty}}\|u_x\|_{L^{\infty}}\|u_{xx}\|^{2n}_{L^{2n}},
\end{split}
\end{equation*}
and
\begin{equation*}
\begin{split}
\left|\int_{\mathbb{R}}u_{xx}^{2n-1}u_x^2u_{xxx}\ dx\right| &=
\left|\dfrac{1}{2n}\int_{\mathbb{R}}u_x^2(u_{xx}^{2n})_x\ dx\right|=\left|-\dfrac{1}{2n}\int_{\mathbb{R}}2u_xu_{xx}^{2n+1}\ dx\right|\\
&\le\dfrac{2}{2n}\|u_x\|_{L^{\infty}}\|u_{xx}\|_{L^{\infty}}\|u_{xx}\|^{2n}_{L^{2n}}.
\end{split}
\end{equation*}
Similarly, one gets
\begin{equation*}
\begin{split}
&\left|\int_{\mathbb{R}}uu_xu_{xx}^{2n}\ dx\right|
\le\|u\|_{L^{\infty}}\|u_x\|_{L^{\infty}}\|u_{xx}\|^{2n}_{L^{2n}},\,
\left|\int_{\mathbb{R}}u_x^3u_{xx}^{2n-1}\
dx\right|\le\|u_x\|^2_{L^{\infty}}\|u_{xx}\|^{2n-1}_{L^{2n}}\|u_x\|_{L^{2n}},\\
&\left|\int_{\mathbb{R}}u^2u_xu_{xx}^{2n-1}\
dx\right|\le\|u\|^2_{L^{\infty}}\|u_{xx}\|^{2n-1}_{L^{2n}}\|u_x\|_{L^{2n}},\,
\left|\int_{\mathbb{R}}u_xu_{xx}^{2n+1}\,dx\right|
\le\|u_x\|_{L^{\infty}}\|u_{xx}\|_{L^{\infty}}\|u_{xx}\|^{2n}_{L^{2n}},
\end{split}
\end{equation*}
and
\begin{equation*}
\begin{split}
 &\left|\int_{\mathbb{R}}u_{xx}^{2n-1}
G\ast\left(\partial_x\left(\frac{2}{3}u^3+uu^2_x\right)+\frac1{3}u^3_x+\gamma
u_x\right)\ d
x\right|\\&\quad\le\|u_{xx}\|^{2n-1}_{L^{2n}}\left(\frac1{3}\|u_x\|_{L^{\infty}}^2\|u_x\|_{L^{2n}}+
\frac{2}{3}\|u\|_{L^{\infty}}^2\|u\|_{L^{2n}}+\|u_x\|_{L^{\infty}}^2\|u\|_{L^{2n}}+|\gamma|\|u_x\|_{L^{2n}}\right).
\end{split}
\end{equation*}
Integrating the resulting inequality over $[0,t]$ and letting $n$ tend to infinity, it
follows that
\begin{align}\label{4.7}
\|u_{xx}(t)\|_{L^{\infty}}\le&\|u_{xx}(0)\|_{L^{\infty}}\nonumber\\&+\int^t_0\|u(\tau)\|_{L^{\infty}}
\left(\frac{2}{3}\|u(\tau)\|_{L^{\infty}}^2+\frac{7}{2}\|u_x(\tau)\|_{L^{\infty}}^2+\frac{5}{2}\|u_{xx}(\tau)\|_{L^{\infty}}^2\right)\
d\tau\nonumber\\&+\int^t_0\|u_x(\tau)\|_{L^{\infty}}(2\|u(\tau)\|_{L^{\infty}}^2+
\|u_x(\tau)\|_{L^{\infty}}^2+3\|u_{xx}(\tau)\|_{L^{\infty}}^2+2|\gamma|)\
d\tau.
\end{align}
If $\gamma\neq0$, let
$$h(t)=2\|u(t)\|_{L^{\infty}}+\|u_{x}(t)\|_{L^{\infty}}+2\|u_{xx}(t)\|_{L^{\infty}}.$$
Combining \eqref{4.3}, \eqref{4.5} and \eqref{4.7}, we deduce that
\begin{equation}\label{4.8}
\|m(t)\|_{L^{\infty}}\le h(t)\le
h(0)+5\int^t_0\left[(h(\tau))^3+2|\gamma|h(\tau)\right]\ d\tau.
\end{equation} Define
$$T=\dfrac1{20|\gamma|}\ln\left(1+\dfrac{2|\gamma|}{(2\|u_0\|_{L^{\infty}}+\|\pa_xu_0\|_{L^{\infty}}
+2\|\pa^2_xu_0\|_{L^{\infty}})^2}\right).$$ By \eqref{4.8}, then for
all $t\le\min\{T,T^{\ast}\}$, one can easily get
$$\|m(t)\|_{L^{\infty}}\le\dfrac{\sqrt{2|\gamma|}h(0)}{\sqrt{(h^2(0)+2|\gamma|)e^{-20|\gamma|t}-h^2(0)}}.$$
By virtue of Theorem \ref{t3.1}, it follows that $T^{\ast}\ge T$.

If $\gamma=0$, let
$$h(t)=3\|u(t)\|_{L^{\infty}}+3\|u_{x}(t)\|_{L^{\infty}}+\|u_{xx}(t)\|_{L^{\infty}}.$$
Combining \eqref{4.3}, \eqref{4.5} and \eqref{4.7}, we obtain that
\begin{equation}\label{4.9}
\|m(t)\|_{L^{\infty}}\le h(t)\le h(0)+\int^t_0(h(\tau))^3\ d\tau.
\end{equation} And let
$$T=\dfrac1{2(3\|u_0\|_{L^{\infty}}+3\|\pa_xu_0\|_{L^{\infty}}
+\|\pa^2_xu_0\|_{L^{\infty}})^2}.$$ Similarly, we obtain that for
all $t\le\min\{T,T^{\ast}\}$,
\begin{equation}\label{4.10}
\|m(t)\|_{L^{\infty}}\le h(t)\le \dfrac{h(0)}{\sqrt{1-8h^2(0)t}}.
\end{equation}
 This completes the proof of Theorem \ref{t4.2}.
\end{proof}

\renewcommand{\theequation}{\thesection.\arabic{equation}}
\setcounter{equation}{0}

\section{Blow-up data for $\gamma=0$}
In this section, we will provide sufficient conditions for the blow-up data to the initial-value problem
 \eqref{e1.1} with $ \gamma = 0. $ The blow-up result is now established in the following.
\begin{thm}\label{t5-1}
Let $ \gamma = 0. $ Suppose $u_0\in H^s\cap L^1 $ with $ s > 5/2$.
Let $T>0$ be the maximal time of existence of the corresponding
solution $m(t,x)$ to \eqref{e1.1} with the initial data  $
m_0(x)=(1-\partial_x^2)u_0$. Assume $m_0(x) \geq 0$ for all $x\in
\mathbb{R}$ and $m_0(x_0)>0$ at some point $x_0 \in \mathbb{R}$.
\begin{itemize}
\item[i)] \,  If
  \begin{equation}\label{blow-4-0}
\partial_x u_0(x_0) <-\|u_0\|_{H^1} \sqrt{\frac{I_0}{m_0(x_0)}} \quad \mbox{with} \quad I_0:=\int_{\mathbb{R}} u_0(x)\,dx,
\end{equation}
then the solution $m(t, x)$ blows up at a time
 \begin{equation*}\label{blow-4-1}
T_0 \leq t^{\ast}:=\frac{- \partial_x
u_0(x_0)}{I_0\|u_0\|_{H^{1}}^2 }-
\sqrt{\left (\frac{\partial_x
u_0(x_0)}{I_0\|u_0\|_{H^{1}}^2 } \right )^2-\frac{1}{I_0\|u_0\|_{H^{1}}^2 {m}_0(x_0)}}.
\end{equation*}
Moreover when $T_0=t^{\ast}$, the following estimate of the blow-up
rate holds
\begin{equation}\label{blow-4-1-a}
\begin{split}
&\liminf_{t\rightarrow T_0^{-}}\left((T_0-t)\,\inf_{x\in
\mathbb{R}}(mu_x)(t, x)\right) \leq - \frac{1}{2}.
\end{split}
\end{equation}

\item[ii)] \,  If
  \begin{equation}\label{blow-4-0-b}
\partial_x u_0(x_0) > - I_0 \quad \mbox{and} \quad \frac{1}{m_{0}(x_0)}-\frac{\partial_x u_0(x_0)}{\sqrt{2}I_0
 \|u_0\|_{H^1}} < \frac{1}{\sqrt{2} \|u_0\|_{H^1}} \ln \bigg( \frac{I_0}{I_0+ \partial_x u_0(x_0) }\bigg),
\end{equation}
then the solution $m(t, x)$ blows up at a time
 \begin{equation*}\label{blow-4-1-b}
T_0 \leq t^{\ast\ast}:=\frac{1}{\sqrt{2} I_0\|u_0\|_{H^1}} \ln \bigg( \frac{I_0}{I_0+ \partial_x u_0(x_0) }\bigg).
\end{equation*}

\item[iii)] \,  If $ \displaystyle
\partial_x u_0(x_0) \le - I_0,
$  then the solution $m(t, x)$ blows up at a time $
T_0 \leq t_1,$
where $t_1$ uniquely solves the equation
 \begin{equation*}
\frac{\sqrt{2}(I_0+ \partial_x u_0(x_0))}{4I_0\|u_0\|_{H^1}}(e^{\sqrt{2}I_0\|u_0\|_{H^1}\,t}-1)- I_0\,t +\frac{1}{m_0(x_0)}=0.
 \end{equation*}
\end{itemize}
\end{thm}

\begin{proof} Denoting $M=m \, u_x$, we first recall from Proposition 5.1 in \cite{gloq} that for all $(t, x) \in [0, T) \times \mathbb{R}$
\begin{equation*}\label{blow-1-0}
 M_{t}+(u^2-u_x^2)M_x=-2M^2-2m (1-\partial_x^2)^{-1}(u_x^2 m)-2m\partial_x(1-\partial_x^2)^{-1}(uu_xm),
\end{equation*}
which along with \eqref{flow-1} implies
 \begin{equation}\label{blow-1-0-a-1}
 \frac{d}{dt}M(t, q(t, x))=\bigg(-2M^2-2m (1-\partial_x^2)^{-1}(u_x^2 m)-2m\partial_x(1-\partial_x^2)^{-1}(uu_xm)\bigg)(t, q(t, x)).
\end{equation}
Since $m_0(x) \geq 0$ for all $x
\in \mathbb{R}$,  Remark \ref{rmk-4-1} implies that
 \begin{equation}\label{blow-2-2}m(t, x) \geq 0,
\end{equation}
for all $t \in [0, T)$, $x \in \mathbb{R}$, and hence
\begin{equation*}\label{blow-2-2-a}\Big(m(1-\partial_x^2)^{-1}(u_x^2 m)\Big)(t, x) \geq 0.
\end{equation*}
On the other hand, for $G=\frac{1}{2}e^{-|x|}$, we have
 \begin{equation*}
 \begin{split}
 \partial_x (1-\partial_x^2)^{-1}(uu_x m)(t, x) &=\partial_x G \ast (uu_x m)(t, x)\\
 &=-\frac{1}{2}\int_{-\infty}^{+\infty}\sign (x-y) e^{-|x-y|}(uu_xm)(t,y)\,dy,
\end{split}
\end{equation*}
which implies
 \begin{equation}\label{blow-1-0-a2}
 \begin{split}
-2m\partial_x(1-\partial_x^2)^{-1}(uu_xm)&=m\int_{-\infty}^{+\infty}\sign(x-y)e^{-|x-y|}(uu_xm)(y) \, dy\\
&\leq m\int_{-\infty}^{+\infty} e^{-|x-y|}(u|u_x|m)(y) \, dy \\&=
2m(1-\partial_x^2)^{-1}(u|u_x|m).
\end{split}
\end{equation}
Therefore, we find from \eqref{blow-1-0-a-1} and \eqref{blow-1-0-a2} that
 \begin{equation}\label{blow-1-0-a-3}
 \frac{d}{dt}M(t, q(t, x))\leq\bigg(-2M^2+2m (1-\partial_x^2)^{-1}((u-|u_x|)|u_x| m)\bigg)(t, q(t, x)).
\end{equation}
Notice that
 \begin{equation*}u(t, x)=(G\ast m)(t, x)=\frac{1}{2}\int_{\mathbb{R}} e^{-|x-y|}m(t,
y)\,dy,
\end{equation*}
we have
 \begin{equation*}\label{blow-2-3}
\begin{split}u(t, x)
 &=\frac{e^{-x}}{2}\int_{-\infty}^{x} e^{y} m(t, y)\, dy
 +\frac{e^{x}}{2}\int_{x}^{+\infty} e^{-y}m(t, y)\,dy,
\\ u_x(t, x)
 &=-\frac{e^{-x}}{2}\int_{-\infty}^{x} e^{y} m(t, y)\, dy
 +\frac{e^{x}}{2}\int_{x}^{+\infty} e^{-y}m(t, y)\,dy,
 \end{split}
\end{equation*}
which along with \eqref{blow-2-2} and Remark \ref{rmk-conservation-1} leads to
\begin{equation*}\label{blow-2-5}
\begin{split}
&0 \leq \int_{x}^{+\infty} e^{x-y}m(t, y)\,dy=u(t, x)+ u_x(t, x) \leq \int_{-\infty}^{+\infty} m(t, y)\,dy=I_0,\\
&0 \leq \int^{x}_{-\infty} e^{y-x}m(t, y)\,dy=u(t, x)-u_x(t, x)\leq \int_{-\infty}^{+\infty} m(t, y)\,dy=I_0.
\end{split}
\end{equation*}
From this, we obtain for all $(t, x) \in [0, T) \times \mathbb{R}$
\begin{equation}\label{blow-2-6}
|u_x(t, x)| \leq u(t, x), \quad u(t, x)-|u_x(t, x)| \leq I_0 \quad \mbox{and} \quad u(t, x) \leq I_0+u_x(t, x).
\end{equation}
Therefore, in view of the Sobolev inequality
$\|u\|_{L^{\infty}(\mathbb{R})} \leq \frac{1}{\sqrt{2}}
\|u\|_{H^1{(\mathbb{R})}}$ and \eqref{blow-2-6}, it follows from \eqref{blow-1-0-a-3} that
 \begin{equation*}\label{blow-1-0-a-4}
  \begin{split}
 \frac{d}{dt}M(t, q(t, x))&\leq -2M^2(t, q(t, x))+2I_0 \|u\|_{L^{\infty}} m(t, q(t, x))[(1-\partial_x^2)^{-1}m](t, q(t, x))\\
 & =-2M^2(t, q(t, x))+2I_0 \|u\|_{L^{\infty}} (mu)(t, q(t, x))\\
 &\leq -2M^2(t, q(t, x))+2I_0 \|u\|_{L^{\infty}}^2 m(t, q(t, x))\\
 &\leq -2M^2(t, q(t, x))+I_0 \|u_0\|_{H^{1}}^2 m(t, q(t, x))
 \end{split}
\end{equation*}
and
 \begin{equation*}\label{blow-1-0-a-5}
  \begin{split}
 \frac{d}{dt}M(t, q(t, x))&\leq  -2M^2(t, q(t, x))+2I_0 \|u\|_{L^{\infty}} (mu)(t, q(t, x))\\
 &\leq -2M^2(t, q(t, x))+2I_0 \|u\|_{L^{\infty}}[m(I_0+u_x)](t, q(t, x))\\
 &\leq -2M^2(t, q(t, x))+\sqrt{2}I_0 \|u_0\|_{H^{1}} M(t, q(t, x))+\sqrt{2}I_0^2 \|u_0\|_{H^{1}} m(t, q(t, x)),
 \end{split}
\end{equation*}
which, in particular, implies
 \begin{equation}\label{blow-1-0-a-4-1}
  \begin{split}
 \frac{d}{dt}M(t, q(t, x_0))\leq -2M^2(t, q(t, x_0))+I_0 \|u_0\|_{H^{1}}^2 m(t, q(t, x_0))
 \end{split}
\end{equation}
and \begin{equation}\label{blow-1-0-a-6}
  \begin{split}
 \frac{d}{dt}M(t, q(t, x_0))\leq -2M^2(t, q(t, x_0))&+\sqrt{2}I_0 \|u_0\|_{H^{1}} M(t, q(t, x_0))\\
 &+\sqrt{2}I_0^2 \|u_0\|_{H^{1}} m(t, q(t, x_0)).
 \end{split}
\end{equation}

Similarly, one can see from the equation in \eqref{e1.1} that
 \begin{equation}\label{blow-4-2-a}
 \begin{split}
\frac{d}{dt}m(t, q(t, x_0))= -2m M(t, q(t, x_0)). \end{split}
\end{equation}
Denote that
$$\overline{M}(t):=2M(t, q(t, x_0)) \quad \mbox{and}\quad
\overline{m}(t):= 2 m(t, q(t, x_0)).$$

We first reformulate \eqref{blow-1-0-a-4-1} and \eqref{blow-4-2-a} as
\begin{equation*}\label{blow-4-3-a-1}
 \frac{d}{dt}\overline{M}(t)\leq -\overline{M}(t)^2+I_0 \|u_0\|_{H^{1}}^2 \overline{m}(t)
\end{equation*}
and
\begin{equation}\label{blow-4-3-a}
\frac{d}{dt}\overline{m}(t)= -\overline{m}(t)\overline{M}(t).
\end{equation}
Combining this with \eqref{blow-2-2}, we deduce that
\begin{equation*}\label{blow-4-4-app-1}
 \begin{split}
\frac{d}{dt}&\left(\frac{1}{\overline{m}(t)^2}\right.\left.\frac{d}{dt}\overline{m}(t)\right)
=\frac{d}{dt}\left(-\frac{1}{\overline{m}(t)}\overline{M}(t)\right)\\
&=\frac{1}{\overline{m}(t)^2}\left(-\overline{m}(t)\frac{d}{dt}\overline{M}(t)
+\overline{M}(t)\frac{d}{dt}\overline{m}(t)\right)\\
&\geq  \frac{1}{\overline{m}(t)^2}
\Big(\overline{m}(t)\big(\overline{M}(t)^2-I_0 \|u_0\|_{H^{1}}^2 \overline{m}(t)\big)
-\overline{m}(t)\overline{M}(t)^2\Big)=-I_0\|u_0\|_{H^{1}}^2.
\end{split}
\end{equation*}

Integrating from $0$ to $t$ leads to
\begin{equation}\label{blow-4-5}
 \begin{split}
\frac{1}{\overline{m}(t)^2}\frac{d}{dt}{\overline{m}(t)} \geq
C_0-I_0\|u_0\|_{H^{1}}^2 t,
\end{split}
\end{equation}
with
\begin{equation*}C_0:=-\frac{\overline{M}(0)}{\overline{m}(0)}
=-(\partial_xu_0)(x_0).
\end{equation*}
Combining this with \eqref{blow-4-3-a} yields
\begin{equation}\label{blow-4-6}
 \begin{split}
\overline{M}(t)
=-\frac{1}{\overline{m}(t)}\frac{d}{dt}{\overline{m}(t)}\leq
-\overline{m}(t)\big(C_0-I_0\|u_0\|_{H^{1}}^2 t\big).
\end{split}
\end{equation}
Integrating \eqref{blow-4-5} again on $[0, t]$ implies
\begin{equation*}\label{blow-4-5-a}
 \begin{split}
\frac{1}{\overline{m}(t)}-\frac{1}{\overline{m}(0)}\leq \frac{1}{2}I_0\|u_0\|_{H^{1}}^2 t^2-C_0
t,
\end{split}
\end{equation*}
and hence
\begin{equation*}\label{blow-4-7}
 \begin{split}
\frac{1}{\overline{m}(t)}&\leq
\frac{1}{2}I_0\|u_0\|_{H^{1}}^2 \left(t^2-\frac{2C_0}{I_0\|u_0\|_{H^{1}}^2 }t+\frac{2}{I_0\|u_0\|_{H^{1}}^2 \overline{m}(0)}\right)\\
&= \frac{1}{2}I_0\|u_0\|_{H^{1}}^2 \left(t^2-\frac{2C_0}{I_0\|u_0\|_{H^{1}}^2 }t+\frac{1}{I_0\|u_0\|_{H^{1}}^2 {m}_0(x_0)}\right).
\end{split}
\end{equation*}
The quadratic equation
\begin{equation*}t^2-\frac{2C_0}{I_0\|u_0\|_{H^{1}}^2 }t+\frac{1}{I_0\|u_0\|_{H^{1}}^2 {m}_0(x_0)}=0
\end{equation*}
has two roots:
\begin{equation*}
 \begin{split}
 t^{\ast}&:=\frac{C_0}{I_0\|u_0\|_{H^{1}}^2 }-
\sqrt{\left (\frac{C_0}{I_0\|u_0\|_{H^{1}}^2 } \right
)^2-\frac{1}{I_0\|u_0\|_{H^{1}}^2 {m}_0(x_0)}}, \qquad \text{and} \\
t_{\ast}&:=\frac{C_0}{I_0\|u_0\|_{H^{1}}^2 }+ \sqrt{\left
(\frac{C_0}{I_0\|u_0\|_{H^{1}}^2 } \right
)^2-\frac{1}{I_0\|u_0\|_{H^{1}}^2 {m}_0(x_0)}}
\end{split}
\end{equation*}
It thus transpires from Assumption \eqref{blow-4-0} that
\begin{equation*}\label{blow-4-8}
\left (\frac{C_0}{I_0\|u_0\|_{H^{1}}^2 } \right )^2>\frac{1}{I_0\|u_0\|_{H^{1}}^2 {m}_0(x_0)}\,, \qquad \hbox{hence} \qquad
0< t^{\ast} < \frac{C_0}{I_0\|u_0\|_{H^{1}}^2 }< t_{\ast}.
\end{equation*}
Thus,
\begin{equation}\label{blow-4-9-0}
 \begin{split}
0 \leq \frac{1}{\overline{m}(t)}\leq \frac{I_0\|u_0\|_{H^{1}}^2}{2}(t-t^{\ast})(t-t_{\ast}).
\end{split}
\end{equation}
It is then adduced from \eqref{blow-4-9-0} that there is a time $ T_0 \in (0, t^{\ast}]$
such that
\begin{equation*}{m}(t) \longrightarrow +\infty, \quad  \mbox{as} \quad t \longrightarrow T_0 \leq t^{\ast},
\end{equation*}
which, by \eqref{blow-4-6}, implies that
\begin{equation*}M(t) \longrightarrow  -\infty, \quad  \mbox{as} \quad t \longrightarrow  T_0 \leq t^{\ast}.
\end{equation*}
Therefore,
\begin{equation*}
\inf_{x\in \mathbb{R}}M(t, x) \leq M(t)\longrightarrow  -\infty, \quad
\mbox{as} \quad t \longrightarrow  T_0 \leq t^{\ast},
\end{equation*}
which, in view of Theorem \ref{t4.2}, implies
that the solution $m(t, x)$ blows up at the time $T_0$.

Having established wave breaking results for \eqref{e1.1} as above,
attention is now given to blow-up rate for the solution. In fact,
owing to \eqref{blow-4-6} and \eqref{blow-4-9-0}, we derive that for
all $0<t<T_0$
\begin{equation*}\label{blow-4-11}
\begin{split}
(T_0-t)\,\inf_{x\in \mathbb{R}}\overline{M}(t, x) &\leq (T_0-t)\, \overline{M}(t) \leq
(T_0-t)\,\overline{m}(t)((\partial_xu_0)(x_0)+I_0\|u_0\|_{H^{1}}^2 t)\\
& \leq
(T_0-t)\,\frac{2}{I_0\|u_0\|_{H^{1}}^2(t-t^{\ast})(t-t_{\ast})}(I_0\|u_0\|_{H^{1}}^2 t+(\partial_xu_0)(x_0))\\
& \leq
\frac{2(T_0-t)}{(t-t^{\ast})(t-t_{\ast})} \left (t+\frac{(\partial_xu_0)(x_0)}{I_0\|u_0\|_{H^{1}}^2} \right),
\end{split}
\end{equation*}
which leads to \eqref{blow-4-1-a} when $T_0=t^{\ast}$. Therefore, we
end the proof of Theorem \ref{t5-1} (i).

On the other hand, we reformulate \eqref{blow-1-0-a-6} as
\begin{equation*}\label{blow-4-3}
 \frac{d}{dt}\overline{M}(t)\leq -\overline{M}(t)^2+\sqrt{2}I_0 \|u_0\|_{H^{1}} \overline{M}(t)+ \sqrt{2} I_0^2 \|u_0\|_{H^{1}} \overline{m}(t)
\end{equation*}
Combining this with \eqref{blow-4-3-a} and \eqref{blow-2-2}, we deduce that
\begin{equation*}\label{blow-4-4}
 \begin{split}
\frac{d}{dt}&\left(-\frac{\overline{M}(t)}{\overline{m}(t)}\right)=\frac{d}{dt}\left(\frac{1}{\overline{m}(t)^2}\right.
\left.\frac{d}{dt}\overline{m}(t)\right)
=\frac{1}{\overline{m}(t)^2}\left(-\overline{m}(t)\frac{d}{dt}\overline{M}(t)
+\overline{M}(t)\frac{d}{dt}\overline{m}(t)\right)\\
&\geq  \frac{1}{\overline{m}(t)^2}
\Bigg[\overline{m}(t)\bigg(\overline{M}(t)^2-\sqrt{2}I_0 \|u_0\|_{H^{1}} \overline{M}(t)- \sqrt{2}I_0^2 \|u_0\|_{H^{1}}
 \overline{m}(t)\bigg)-\overline{m}(t)\overline{M}(t)^2\Bigg]\\
&=-\frac{\overline{M}(t)}{\overline{m}(t)}\sqrt{2}I_0 \|u_0\|_{H^{1}}- \sqrt{2} I_0^2 \|u_0\|_{H^{1}},
\end{split}
\end{equation*}
which gives rise to
\begin{equation}\label{blow-4-8}
 \begin{split}
&\frac{d}{dt}\left(-\frac{\overline{M}(t)}{\overline{m}(t)}e^{-C_1\,t}\right)
\geq -C_2e^{-C_1\,t}
\end{split}
\end{equation}
with
\begin{equation*}
  C_1:= \sqrt{2}I_0 \|u_0\|_{H^{1}}, \quad C_2:= \sqrt{2} I_0^2 \|u_0\|_{H^{1}}.
\end{equation*}
Integrating \eqref{blow-4-8} from $0$ to $t$ leads to
\begin{equation}\label{blow-4-9}
 \begin{split}
-\frac{\overline{M}(t)}{\overline{m}(t)}e^{-C_1\,t} \geq
\frac{C_2}{C_1}e^{-C_1\,t}+\frac{C_0C_1-C_2}{C_1}
\end{split}
\end{equation}
with
\begin{equation*}C_0=-\frac{\overline{M}(0)}{\overline{m}(0)}
=-(\partial_xu_0)(x_0),
\end{equation*}
which implies
\begin{equation}\label{blow-4-9-a-1}
 \begin{split}
\overline{M}(t) \leq -\overline{m}(t)\bigg(\frac{C_2}{C_1}+\frac{C_0C_1-C_2}{C_1}e^{C_1\,t}\bigg)
\end{split}
\end{equation}
and
\begin{equation}\label{blow-4-10}
 \begin{split}
-\frac{d}{dt}\bigg(\frac{1}{\overline{m}(t)}\bigg)\geq
\frac{C_2}{C_1}+\frac{C_0C_1-C_2}{C_1}e^{C_1\,t}.
\end{split}
\end{equation}
Integrating \eqref{blow-4-10} again on $[0, t]$ implies
\begin{equation*}\label{blow-4-11}
 \begin{split}
\frac{1}{\overline{m}(t)}-\frac{1}{\overline{m}(0)}\leq \frac{C_2-C_0C_1}{C_1^2}(e^{C_1\,t}-1)-\frac{C_2}{C_1}\,t,
\end{split}
\end{equation*}
and hence
\begin{equation*}\label{blow-4-12}
 \begin{split}
0 \leq \frac{1}{\overline{m}(t)}\leq \frac{C_2-C_0C_1}{C_1^2}(e^{C_1\,t}-1)-\frac{C_2}{C_1}\,t +\frac{1}{\overline{m}(0)}=: f(t).
\end{split}
\end{equation*}
Notice that $f(0)=\frac{1}{\overline{m}(0)} >0$ and the assumption \eqref{blow-4-0-b} implies the equation $\frac{d}{dt}f(t)=0$ has only one root
\begin{equation*}
 \begin{split}t^{\ast\ast}:=\frac{1}{C_1} \ln \left (\frac{C_2}{C_2-C_0C_1} \right )
\end{split}
\end{equation*}
and then
\begin{equation*}\label{blow-4-13}
f(t^{\ast\ast})= \frac{C_0}{C_1}+\frac{1}{m_0(x_0)}-\frac{C_2}{C_1^2} \ln \left (\frac{C_2}{C_2-C_0C_1} \right ) < 0.
\end{equation*}
From this, we may find a time $0<T_0 \leq t^{\ast\ast}$
such that
\begin{equation*}{m}(t) \longrightarrow +\infty, \quad  \mbox{as} \quad t \longrightarrow T_0 \leq t^{\ast\ast},
\end{equation*}
which, by \eqref{blow-4-9-a-1}, implies that
\begin{equation*}M(t) \longrightarrow  -\infty, \quad  \mbox{as} \quad t \longrightarrow  T_0 \leq t^{\ast\ast}.
\end{equation*}
Therefore,
\begin{equation*}
\inf_{x\in \mathbb{R}}M(t, x) \leq M(t)\longrightarrow  -\infty, \quad
\mbox{as} \quad t \longrightarrow  T_0 \leq t^{\ast\ast},
\end{equation*}
which, in view of Theorem \ref{t4.2}, implies
that the solution $m(t, x)$ blows up at the time $T_0$. This
ends the proof of Theorem \ref{t5-1}(ii). Similarly, we may prove  Theorem \ref{t5-1} (iii). This completes the proof of  Theorem \ref{t5-1}.
\end{proof}

\begin{rmk}\label{rmk-5-1}
(1) Compared to the blow-up result of Theorem 5.2 in \cite{gloq}
where $\partial_x u_0(x_0)
<-\sqrt{\frac{\sqrt{2}\|u_0\|_{H^1}^3}{m_0(x_0)}}$, Theorem
\ref{t5-1}(i) looks better at least in some sense as the following
example: taking the initial data $u_0(x)=e^{-x^2}$, then
\begin{equation*}I_0=\int_{-\infty}^{\infty}e^{-x^2}\, dx=\sqrt{\pi},\quad \|u_0\|_{H^1}^2=\int_{-\infty}^{\infty}e^{-2x^2}\, dx
+4\int_{-\infty}^{\infty}x^2e^{-2x^2}\, dx=\sqrt{2\pi},
\end{equation*}
which implies
$$\|u_0\|_{H^1} \sqrt{\frac{I_0}{m_0(x_0)}} < \sqrt{\frac{\sqrt{2}\|u_0\|_{H^1}^3}{m_0(x_0)}}.$$

(2) In view of  Theorem \ref{t5-1} (iii), if $ \partial_x
u_0(x_0)\le - I_0$, then the function
 \begin{equation*}
F(t):=\frac{\sqrt{2}(I_0+ \partial_x u_0(x_0))}{4I_0\|u_0\|_{H^1}}(e^{\sqrt{2}I_0\|u_0\|_{H^1}\,t}-1)- I_0\,t +\frac{1}{m_0(x_0)}=0.
 \end{equation*}
 has only one root on $[0, +\infty)$. In fact, we need only consider the case $ \partial_x u_0(x_0) < - I_0$, that is,
  $I_0+ \partial_x u_0(x_0)<0$. Notice that $F(0)=\frac{1}{m_0(x_0)}>0$, $F(+\infty)<0$, and $\frac{d}{dt}F(t)<0$ for all
   $t \in [0, +\infty)$.  We then deduce from the Intermediate Value Theorem that $F(t)=0$ has only one root on $[0, +\infty)$.
\end{rmk}
\renewcommand{\theequation}{\thesection.\arabic{equation}}
\setcounter{equation}{0}

\section{Non-existence of smooth traveling waves for $\gamma=0$}

In this section, we prove that the equation in  (\ref{e1.1}) does
not have nontrivial smooth traveling waves.

\begin{thm}\label{t6.1}
There is no nontrivial smooth traveling wave solution
$u(t,x)=\phi(x-ct)$, $c\in {\Bbb R}$ of the Cauchy problem
(\ref{e1.1}) with $\gamma=0$ in $C([0,\infty); H^3({\Bbb R}))\cap
C^1([0,\infty); H^2({\Bbb R}))$.
\end{thm}
\begin{proof}  We use a contradiction argument. Assume that $\phi\in
H^3$ is a strong solution of the Cauchy problem (\ref{e1.1}). Then
we have
\begin{equation*}
c(\phi-\phi'')'=((\phi^2-\phi_x^2)(\phi-\phi''))'\;\; {\rm in}\;\;
L^2({\Bbb R}).
\end{equation*}
Since $\phi\in H^3({\Bbb R})\subset C_0^2({\Bbb R})$, we find that
\begin{equation}\label{e6.2}
c(\phi-\phi'')=(\phi^2-\phi_x^2)(\phi-\phi'')\;\; {\rm in}\;\;
H^1({\Bbb R}).
\end{equation}
Note that $\phi \not\equiv 0$ and $\phi,\phi',\phi''\rightarrow 0$
as $|x|\rightarrow \infty$, it implies that $\phi-\phi''\neq 0$.
Otherwise, $\phi=c_1e^x+c_2 e^{-x}$, which gives $\phi\equiv 0$,
$x\in {\Bbb R}$ since $\phi\rightarrow 0$ as $|x|\rightarrow
\infty$. It then follows from (\ref{e6.2}) that
\begin{equation}\label{e6.3}
\phi^2-\phi'^2=c.
\end{equation}
Let $|x|\rightarrow \infty$. Then $\phi,\phi'\rightarrow 0$. It
yields from (\ref{e6.3}) that $c=0$. Hence we deduce from
(\ref{e6.3}) that
\begin{equation*}
\phi^2-\phi'^2=0,
\end{equation*}
which implies that either $\phi=c_1 e^x$ or $\phi=c_2 e^{-x}$. This then leads to a contradiction that $ \phi \equiv 0, $
  since  $ \phi \to 0, $ as $ |x| \to \infty.$  This completes the proof of the theorem.
\end{proof}

\renewcommand{\theequation}{\thesection.\arabic{equation}}
\setcounter{equation}{0}
\section{The persistence property of solutions}
Attention in this section is now turned to persistence properties of solutions.
In order to derive the persistence properties of solutions, we now  reformulate the equation in (\ref{e1.1}) to the form suitable
for our purpose.

In terms of \eqref{4.2}, the equation in (\ref{e1.1}) can be
rewritten as
\begin{equation}\label{2.1}
u_t+u^2u_x+\partial_xG*\left(\frac{2}{3}u^3+uu^2_x+\gamma
u\right)+\frac1{3}\partial^2_xG*u^3_x=0.
\end{equation}
Our persistence property result for equation (\ref{2.1}) is as
follows:
\begin{thm}\label{t7.1}
Assume that for some $T>0$ and $s>5/2$, $u(t,x)\in C([0,T);H^s)\cap
C^1([0,T);H^{s-1})$, is a strong solution of the initial value
problem associated to equation (\ref{2.1}). If $u_0(x)=u(0,x)$
satisfies that for some $\theta>0$
$$|u_0(x)|,\,|\partial_x u_0(x)|\sim O(\me^{-\theta x})\quad as \quad
x\uparrow\infty.$$ Then
$$|u(t,x)|,\,|\partial_x u(t,x)|\sim O(\me^{-\theta x})\quad as \quad
x\uparrow\infty,$$ uniformly in the time interval $[0,T]$.
\end{thm}
\begin{proof}
For simplicity, we introduce the following notations:
$$H_1(u)=\frac{2}{3}u^3+uu^2_x+\gamma
u,\  H_2(u)=\frac1{3}u^3_x,$$
$$M=\underset{t\in[0,T]}{\sup}\|u(t)\|_{H^s},$$
and
$$\varphi_N(x)=\left\{
\begin{array}{ll}
1, &\quad x\le0,\\
\me^{\theta x}, &\quad x\in(0,N),\\
\me^{\theta N},&\quad x\ge N.\end{array}\right.$$
 For all $N>0$ we have $0\le\varphi'_N(x)\le\theta\varphi_N(x).$
Differentiating (\ref{2.1}) in the $x$-variable yields
\begin{equation}\label{2.3}
u_{xt}+u^2u_{xx}+2uu^2_x+\partial^2_xG*H_1(u)+\partial_xG*H_2(u)-u^2_xu_{xx}=0.
\end{equation}
 Multiplying (\ref{2.1}) and
(\ref{2.3}) by $\varphi_N(x)$ we get
\begin{equation}\label{2.4}
\partial_t(u\varphi_N)+(u\varphi_N)uu_x+\varphi_N\partial_xG*H_1(u)+\varphi_N\partial^2_x
G\ast H_2(u)=0.
\end{equation}
\begin{align}\label{2.5}
&\partial_t(u_x\varphi_N)+2(u_x\varphi_N)uu_x+u^2u_{xx}\varphi_N-u_x^2u_{xx}\varphi_N\nonumber\\
&\quad+\varphi_N \partial_x^2G\ast H_1(u)+\varphi_N\partial_x G\ast
H_2(u)=0.
\end{align}
Again, multiplying (\ref{2.4}) by $(u\varphi_N)^{2p-1}$ and
(\ref{2.5}) by $(u_x\varphi_N)^{2p-1}(p\in\mathbb{Z}^+)$,
respectively, then integrating the results in $x$-variable we obtain
\begin{align*}
&\int_{\mathbb{R}}(u\varphi_N)^{2p-1}\partial_t(u\varphi_N) \ d x=
\dfrac{1}{2p}\dfrac{\ d}{\ d t}\|u\varphi_N\|^{2p}_{L^{2p}}
=\|u\varphi_N\|^{2p-1}_{L^{2p}}\frac{\ d}{\ d t}\|u\varphi_N\|_{L^{2p}},\\
&\left|\int_{\mathbb{R}}(u\varphi_N)^{2p}uu_x\ d x\right|\le
\frac1{2}(\|u(t)\|^2_{L^{\infty}}+\|u_x(t)\|^2_{L^{\infty}})\|u\varphi_N\|^{2p}_{L^{2p}},\\
&\left|\int_{\mathbb{R}}(u\varphi_N)^{2p-1}\varphi_N\partial_x G\ast
H_1(u)\ d
x\right|\le\|u\varphi_N\|^{2p-1}_{L^{2p}}\|\varphi_N\partial_x G\ast
H_1(u)\|_{L^{2p}},
\end{align*}
and
\begin{align*}
&\left|\int_{\mathbb{R}}(u\varphi_N)^{2p-1}\varphi_N
\partial_x^2G\ast H_2(u)\ d
x\right|\le\|u\varphi_N\|^{2p-1}_{L^{2p}}\|\varphi_N
\partial_x^2G\ast H_2(u)\|_{L^{2p}}.
\end{align*}
From (\ref{2.4}) we have
\begin{align}\label{2.6}
\dfrac{\ d}{\ d
t}\|u\varphi_N\|_{L^{2p}}\le&\frac1{2}(\|u(t)\|^2_{L^{\infty}}+\|u_x(t)\|^2_{L^{\infty}})\|u\varphi_n\|_{L^{2p}}\nonumber\\&+
\|\varphi_N\partial_x G\ast
H_1(u)(t)\|_{L^{2p}}+\|\varphi_N\partial_x^2 G\ast
H_2(u)(t)\|_{L^{2p}}.
\end{align}
Similarly,
\begin{align*}
&\int_{\mathbb{R}}(u_x\varphi_N)^{2p-1}\partial_t(u_x\varphi_N) \ d
x= \dfrac{1}{2p}\dfrac{\ d}{\ d t}\|u_x\varphi_N\|^{2p}_{L^{2p}}
=\|u_x\varphi_N\|^{2p-1}_{L^{2p}}\frac{\ d}{\ d t}\|u_x\varphi_N\|_{L^{2p}},\\
&\left|\int_{\mathbb{R}}(u_x\varphi_N)^{2p-1}(u_x\varphi_N)2uu_x\,dx
\right| \le
(\|u(t)\|^2_{L^{\infty}}+\|u_x(t)\|^2_{L^{\infty}})\|u_x\varphi_N\|^{2p}_{L^{2p}},
\end{align*}
and
\begin{align*}
&\left|\int_{\mathbb{R}}(u_x\varphi_N)^{2p-1}(u_{xx}\varphi_N)u^2\,dx
\right|=\left|\int_{\mathbb{R}}u^2(u_x\varphi_N)^{2p-1}[\partial_x(u_x\varphi_N)-u_x\varphi'_N]\
dx\right|\\&\quad=\left|\frac1{2p}\int_{\mathbb{R}}u^2\partial_x(u_x\varphi_N)^{2p}\
d x-\int_{\mathbb{R}}u^2(u_x\varphi_N)^{2p-1}u_x\varphi'_N\ d
x\right|\\&\quad=\left|\frac1{2p}\int_{\mathbb{R}}2uu_x(u_x\varphi_N)^{2p}\
d x+\int_{\mathbb{R}}u^2(u_x\varphi_N)^{2p-1}u_x\varphi'_N\ d
x\right|\\&\quad\le
\left[\left(\frac1{2p}+\theta\right)\|u(t)\|^2_{L^{\infty}}+\frac1{2p}\|u_x(t)\|^2_{L^{\infty}}\right]\|u_x\varphi_N\|^{2p}_{L^{2p}}.
\end{align*}
Since
\begin{align*}
\int_{\mathbb{R}}(u_x\varphi_N)^{2p-1}\varphi_Nu_x^2u_{xx}\ d
x&=\frac1{2}\int_{\mathbb{R}}(u_x\varphi_N)^{2p}\ d
(u_x)^2\\&=-\frac1{2}\int_{\mathbb{R}}2p(u_x\varphi_N)^{2p-1}u_x^2(\varphi_Nu_{xx}+u_x\varphi'_N)\
d x,
\end{align*}
we deduce  that
\begin{align*}
\left|\int_{\mathbb{R}}(u_x\varphi_N)^{2p-1}\varphi_Nu_x^2u_{xx}\ d
x\right|&=\frac{p}{p+1}\left|\int_{\mathbb{R}}u_x^2(u_x\varphi_N)^{2p-1}u_x\varphi'_N\
d
x\right|\\&\leq\frac{p\theta}{p+1}\left|\int_{\mathbb{R}}u^2_x(u_x\varphi_N)^{2p}\
d x\right|\le
\frac{p\theta}{p+1}\|u_x(t)\|^2_{L^{\infty}}\|u_x\varphi_N\|^{2p}_{L^{2p}}.
\end{align*}
For the rest terms of (\ref{2.5}), similarly we have
\begin{align*}
&\left|\int_{\mathbb{R}}(u_x\varphi_N)^{2p-1}\varphi_N
\partial_x^2G\ast H_1(u)\ d
x\right|\le\|u_x\varphi_N\|^{2p-1}_{L^{2p}}\|\varphi_N
\partial_x^2G\ast
H_1(u)\|_{L^{2p}},\\
&\left|\int_{\mathbb{R}}(u_x\varphi_N)^{2p-1}\varphi_N
\partial_xG\ast H_2(u)\ d
x\right|\le\|u_x\varphi_N\|^{2p-1}_{L^{2p}}\|\varphi_N
\partial_xG\ast H_2(u)\|_{L^{2p}}.
\end{align*}
It thus follows from (\ref{2.5}) that
\begin{align}\label{2.7}
\dfrac{\ d}{\ d
t}\|u_x\varphi_N\|_{L^{2p}}\le&(3+\theta)(\|u(t)\|^2_{L^{\infty}}+\|u_x(t)\|^2_{L^{\infty}})
\|u_x\varphi_N\|_{L^{2p}}\nonumber\\&+ \|\varphi_N \partial_x^2G\ast
H_1(u)(t)\|_{L^{2p}}+\|\varphi_N
\partial_x G\ast H_2(u)(t)\|_{L^{2p}}.
\end{align}
For any $n>\frac1{2}$, the Sobolev embedding theorem
$H^n\hookrightarrow L^{\infty}$ holds. Thus
$\|u(t)\|_{L^{\infty}}\le\|u(t)\|_{H^s}\le M$ and
$\|u_x(t)\|_{L^{\infty}}\le\|u(t)\|_{H^s}\le M$ hold for $s>5/2$. By
Gronwall's inequality, from (\ref{2.6}) and (\ref{2.7}) we arrive at
\begin{align}\label{2.8}
\|u\varphi_N\|_{L^{2p}}\le&\Big(\|u(0)\varphi_N(0)\|_{L^{2p}}+\int^t_0(\|\varphi_N\partial_x
G\ast H_1(u)\|_{L^{2p}}\nonumber\\&+\|\varphi_N \partial_x^2G\ast
H_2(u)\|_{L^{2p}})\ d\tau\Big)\me^{M^2t},
\end{align}
and
\begin{align}\label{2.9}
\|u_x\varphi_N\|_{L^{2p}}\le&[\|u_x(0)\varphi_N(0)\|_{L^{2p}}\nonumber\\&+\int^t_0(\|\varphi_N\partial^2_x
G\ast H_1(u)(t)\|_{L^{2p}}+\|\varphi_N\partial_x G\ast
H_2(u)(t)\|_{L^{2p}})\ d\tau]\me^{2(3+\theta)M^2t}.
\end{align}
Since $f\in L^2\cap L^{\infty}$ implies
$\underset{q\uparrow\infty}{\lim}\|f\|_{L^q}=\|f\|_{L^{\infty}},$
taking the limits $p\uparrow\infty$ in (\ref{2.8}) and (\ref{2.9})
 and summing the results up we
deduce that
\begin{align}\label{2.10}
\|u\varphi_N&\|_{L^{\infty}}+\|u_x\varphi_N\|_{L^{\infty}}\nonumber\\
\le&\me^{2(3+\theta)M^2t}(\|u(0)\varphi_N\|_{L^{\infty}}+\|u_x(0)\varphi_N\|_{L^{\infty}})\nonumber\\
& +\me^{2(3+\theta)M^2t}\int^t_0(\|\varphi_N\partial_x G\ast
H_1(u)(t)\|_{L^{\infty}}+\|\varphi_N\partial^2_x G\ast
H_1(u)(t)\|_{L^{\infty}}\nonumber\\
& +\|\varphi_N\partial_x G\ast
H_2(u)(t)\|_{L^{\infty}}+\|\varphi_N\partial^2_x G\ast
H_2(u)(t)\|_{L^{\infty}})\ d\tau.
\end{align}

A direct calculation shows that there exists $c>0$, depending on
$\theta>0$ such that for any $N\in\mathbb{Z}^+$
$$\varphi_N(x)\int_{\mathbb{R}}\me^{-|x-y|}
\frac1{\varphi_N(y)}\ d y\le c.$$ Thus, for the terms related to
$H_1(u)$ and $H_2(u)$ we obtain
\begin{align*}
|\varphi_N\partial_x G\ast
H_1(u)|&=\left|\frac1{2}\varphi_N(x)\int_{\mathbb{R}}\sgn(x-y)
\me^{-|x-y|}
H_1(u)\ d y\right|\\
&\le\frac1{2}\varphi_N(x)\left|\int_{\mathbb{R}}\me^{-|x-y
|}\frac1{\varphi_N
(y)}\varphi_N (y) \left(\frac{2}{3}u^3+uu_x^2+\gamma u\right)\ d y\right|\\
&\le c_1(\|u\|^2_{L^{\infty}}+\|
u_x\|^2_{L^{\infty}}+\gamma)\|u\varphi_N\|_{L^{\infty}}.
\end{align*}
Note that $\partial_x G\in L^1,H_1(u),H_2(u)\in L^1\cap
L^{\infty}\quad\text{and}\quad\partial^2_xG*f=(G-\delta)*f$ for any
$f$, one finds that
\begin{align*}
|\varphi_N\partial^2_x G\ast
H_1(u)|&=\left|\frac1{2}\varphi_N(x)\int_{\mathbb{R}}\me^{-|x-y|}
H_1(u)\ d y-\varphi_N(x)\int_{\mathbb{R}}\delta(x-y)H_1(u)\ d y\right|\\
&\le\frac1{2}\varphi_N(x)\left|\int_{\mathbb{R}}\me^{-|x-y|}
\frac1{\varphi_N (y)}\varphi_N (y)
\left(\frac{2}{3}u^3+uu_x^2+\gamma u\right)\ d y\right|
\\&\quad+\left|\varphi_N(x)
\left(\frac{2}{3}u^3+uu_x^2+\gamma u\right)(x)\right|\\
&\le c_2(\|u\|^2_{L^{\infty}}+\|
u_x\|^2_{L^{\infty}}+\gamma)\|u\varphi_N\|_{L^{\infty}}.
\end{align*}
Similarly, we have
\begin{align*}
|\varphi_N\partial_x G\ast
H_2(u)|&=\left|\frac1{2}\varphi_N(x)\int_{\mathbb{R}}\sgn(x-y)
\me^{-|x-y|}
H_2(u)\ d y\right|\\
&\le\frac1{2}\varphi_N(x)\left|\int_{\mathbb{R}}\me^{-|x-y
|}\frac1{\varphi_N
(y)}\varphi_N (y) \frac1{3}u_x^3\ d y\right|\\
&\le c_3\| u_x\|^2_{L^{\infty}}\|u_x\varphi_N\|_{L^{\infty}},
\end{align*}
and
\begin{align*}
|\varphi_N\partial^2_x G\ast
H_2(u)|&=\left|\frac1{2}\varphi_N(x)\int_{\mathbb{R}}\me^{-|x-y|}
H_2(u)\ d y-\varphi_N(x)\int_{\mathbb{R}}\delta(x-y)H_2(u)\ d y\right|\\
&\le\frac1{2}\varphi_N(x)\left|\int_{\mathbb{R}}\me^{-|x-y|}
\frac1{\varphi_N (y)}\varphi_N (y) \frac1{3}u_x^3(y)\ d y\right|
+\frac1{3}|\varphi_N(x)
u_x^3(x)|\\
&\le c_4\| u_x\|^2_{L^{\infty}}\|u_x\varphi_N\|_{L^{\infty}}.
\end{align*}
 Thus, inserting the above estimates into
(\ref{2.10}),
 it follows that there exists a constant $C=C(\theta,M,T)>0$
such that\\
\begin{align*}
\|\varphi_N u\|_{L^{\infty}}+\|\varphi_N u_x\|_{L^{\infty}}
 \le &C(\|\varphi_N u_0\|_{L^{\infty}}+\|\varphi_N
u_x(0)\|_{L^{\infty}})\\&+C\int^t_0 (\|u\|^2_{L^{\infty}}+\|
u_x\|^2_{L^{\infty}}+\gamma )(\|\varphi_N
u\|_{L^{\infty}}+\|\varphi_N u_x\|_{L^{\infty}})\ d\tau.
\end{align*}
Hence, for any $N\in\mathbb{Z}^+$ and any $t\in[0,T]$ we have
\begin{align*}
\|\varphi_N u\|_{L^{\infty}}+\|\varphi_N u_x\|_{L^{\infty}}&\le
C(\|\varphi_N u_0\|_{L^{\infty}}+\|\varphi_N
u_x(0)\|_{L^{\infty}})\\ & \le
 C(\|\me^{\theta x} u_0\|_{L^{\infty}}+\|\me^{\theta x}
u_x(0)\|_{L^{\infty}}).
\end{align*}
Finally, taking the limit as $N$ goes to infinity, we get
$$\underset{t\in[0,T]}{\sup}(\|\me^{\theta x} u\|_{L^{\infty}}+\|\me^{\theta x}
u_x\|_{L^{\infty}})\le C(\|\me^{\theta x}
u_0\|_{L^{\infty}}+\|\me^{\theta x} u_x(0)\|_{L^{\infty}}).$$
Therefore
$$|u(t, x)|,\,|\partial_x u(t,x)|\sim O(\me^{-\theta x})\quad as \quad
x\uparrow\infty,$$ uniformly in the time interval $[0,T]$.
\end{proof}

\begin{thm}\label{t7.2}
Assume that for some $T>0$ and $s>5/2$, $u(x,t)\in C([0,T);H^s) $,
is a strong solution of the initial value problem associated to
equation (\ref{2.1}) with $\gamma=0$. If $u_0(x)=u(x,0)$ satisfies
that for some $\beta\in(\frac1{3},1)$
$$|u_0(x)|\sim O(\me^{-x}),\,|\partial_x u_0(x)|\sim O(\me^{-\beta x})\quad as \quad
x\uparrow\infty,$$  then
$$|u(t, x)|\sim O(\me^{-x})\quad as \quad
x\uparrow\infty,$$ uniformly in the time interval $[0,T]$.
\end{thm}
\begin{proof}
For any $t_1\in(0,T]$, integrating equation (\ref{2.1}) over the
time interval $[0,t_1]$ we obtain
\begin{align*}
u(t_1, x)&-u(0, x)+\int^{t_1}_0u^2u_x(\tau, x)\
d\tau-\frac1{3}\int^{t_1}_0u_x^3(\tau, x)\
d\tau\\&+\int^{t_1}_0\partial_x G\ast
\left(\frac{2}{3}u^3+uu^2_x\right)(\tau, x)\ d\tau+\int^{t_1}_0 G\ast
\frac1{3}u_x^3(\tau, x) d\tau=0.
\end{align*}
By the hypothesis of Theorem \ref{t7.2} we have
$$u(0, x)\sim O(\me^{-x})\quad as \quad x\uparrow\infty.$$
From Theorem \ref{t7.1}, it follows that
$$\int^{t_1}_0u^2u_x(\tau, x)\ d x\sim O(\me^{-(2+\beta) x})\sim o(\me^{-x})\quad as \quad x\uparrow\infty,$$
and
$$\int^{t_1}_0u^3_x(\tau, x)\ d x\sim O(\me^{-3\beta x})\sim o(\me^{-x})\quad as \quad x\uparrow\infty.$$
For the rest two integrals, we deduce that
$$
 \int^{t_1}_0\partial_x
G\ast\left(\frac{2}{3}u^3+uu^2_x\right)(\tau, x)\ d\tau+\int^{t_1}_0 G\ast
\frac1{3}u_x^3(x,\tau) d\tau=\partial_x G\ast H_1(x)+G\ast H_2(x).
$$
Then in view of the definition of Green function $G(x)$, one finds
that
\begin{align*}
\partial_x &G\ast H_1(x)+G\ast H_2(x)\\&=-\frac1{2}\int_{\mathbb{R}}\me^{-|x-y|}\sgn(x-y)
H_1(y)\ d y+\frac1{2}\int_{\mathbb{R}}\me^{-|x-y|}H_2(y)\ d y\\
&=-\frac1{2}\me^{-x}\int^x_{-\infty}\me^y[H_1(y)-H_2(y)]\ d
y+\frac1{2}\me^x\int_x^{\infty}\me^{-y}[H_1(y)+H_2(y)]\ d y.
\end{align*}
Note that
$$H_1(y)\pm H_2(y)\sim O(\me^{-3\beta y})\sim o(\me^{-y}),$$
so we deduce that
$$-\me^{-x}\int^x_{-\infty}\me^y
[H_1(y)-H_2(y)]\ d y\le -C\me^{-x},\quad for \quad x\gg1,$$ and
$$\me^x\int_x^{\infty}\me^{-y}
[H_1(y)+H_2(y)]\ d y\sim o(1)\me^x\int_x^{\infty}\me^{-2y}\ dy\sim
o(\me^{-x}).$$ Therefore,
$$  \int^{t_1}_0\partial_x
G\ast\left(\frac{2}{3}u^3+uu^2_x\right)(\tau, x)\ d\tau+\int^{t_1}_0 G\ast
\frac1{3}u_x^3(\tau, x)\, d\tau \sim O(\me^{-x}).$$ At last, we obtain
$$|u(t_1, x)|\sim O(\me^{-x})\quad as \quad
x\uparrow\infty,$$
 uniformly in the time interval $[0,T]$. Thus considering the arbitrary choose of the time $t_1$,
 this completes the proof of Theorem \ref{t7.2}.
\end{proof}

\renewcommand{\theequation}{\thesection.\arabic{equation}}
\setcounter{equation}{0}

\vskip 0.1cm

\noindent {\bf Acknowledgements.} The work of Fu is partially supported by the NSF-China grant-11001219
and the SPED grant SRP-2010JK860. The work of Gui is partially supported by the NSF of China under the grants
 11001111 and 11171241. The work of Liu is partially
supported by the NSF grants DMS-0906099 and DMS-1207840 and the NHARP
grant-003599-0001-2009. The work of Qu is supported in part by the
NSF-China for Distinguished Young Scholars grant-10925104.

\end{document}